\documentclass[a4paper,11pt]{article}
\usepackage{amssymb,amsmath,amsthm,amscd}
\usepackage[mathscr]{eucal}
\usepackage{fullpage}
\newtheorem{thm}{Theorem}[section]
\newtheorem{lem}[thm]{Lemma}
\newtheorem{cor}[thm]{Corollary}
\newtheorem{prop}[thm]{Proposition}

\theoremstyle{remark}
\newtheorem{rem}[thm]{Remark}
\newtheorem{exmp}[thm]{Example}

\title{On centralizers of parabolic subgroups in Coxeter groups}
\author{Koji Nuida}
\date{}
\begin{document}
\maketitle
\begin{abstract}
Let $W$ be an arbitrary Coxeter group, possibly of infinite rank.
We describe a decomposition of the centralizer $Z_W(W_I)$ of an arbitrary parabolic subgroup $W_I$ into the center of $W_I$, a Coxeter group and a subgroup defined by a $2$-cell complex.
Only information about finite parabolic subgroups is required in an explicit computation.
By using our description of $Z_W(W_I)$, we will be able to reveal a further strong property of the action of the third factor on the second factor, in particular on the finite irreducible components of the second factor.
\end{abstract}
%%%
\section{Introduction}
\label{sec:intro}

\subsection{Background and summary of this work}
\label{sec:intro_background}
A group $W$ is called a Coxeter group if it has a generating set $S$ such that $W$ admits a presentation of the following form:
\begin{displaymath}
W=\langle S \mid (st)^{m_{s,t}}=1 \mbox{ for all } s,t \in S \mbox{ such that } m_{s,t}<\infty \rangle \enspace,
\end{displaymath}
where the data $m_{s,t} \in \{1,2,\dots\} \cup \{\infty\}$ are symmetric in $s,t \in S$ and we have $m_{s,t}=1$ if and only if $s=t$.
The pair $(W,S)$ of a group $W$ and such a generating set $S$ is called a Coxeter system.
Coxeter groups and some associated objects, such as root systems, appear frequently in various topics of mathematics such as geometry, representation theory, combinatorics, etc.
In connection with such related topics, algebraic or combinatorial properties of Coxeter systems and those associated objects have been well studied, forming a long history and establishing many beautiful theories (see e.g., \cite{Hum} and references therein).
However, group structures of Coxeter groups themselves, especially those of infinite Coxeter groups, have not been studied so well; for example, even a fundamental property that an infinite irreducible Coxeter group is always directly indecomposable as an abstract group had not been known until very recent works \cite{Nui_indec,Par_irred}.

The object of this paper is the centralizer $Z_W(W_I)$ of any parabolic subgroup, which is the subgroup $W_I = \langle I \rangle$ generated by a subset $I$ of $S$, in an arbitrary Coxeter group $W$.
We emphasize that we do not place any restriction on $W$ or $W_I$ such as finiteness of cardinalities or finiteness of ranks.
Analogously to other classes of groups such as finite simple groups, it is natural to hope that centralizers of some typical subgroups of Coxeter groups are useful to investigate the structural properties of Coxeter groups as abstract groups.
Indeed, the results of this paper and their enhancements, the latter being planned to be included in a forthcoming paper by the author, are applied in an indispensable manner to the author's recent study on the isomorphism problem in Coxeter groups (see \cite{Nui_ref}).
See Section \ref{sec:intro_related} for other applications.

We give a brief summary of the results of this paper.
We present a decomposition of the centralizer $Z_W(W_I)$ in the following manner:
\begin{displaymath}
Z_W(W_I)=Z(W_I) \times (W^{\perp I} \rtimes B_I) \enspace.
\end{displaymath}
Here the first factor $Z(W_I)$ is the center of $W_I$ whose structure is well known; the second factor $W^{\perp I}$ is the subgroup generated by the reflections along the roots orthogonal to all roots associated to the subsystem $(W_I,I)$ of $(W,S)$; and the third factor $B_I$ is a certain subgroup.
For the second factor $W^{\perp I}$, an existing general theorem implies that such a reflection subgroup $W^{\perp I}$ forms a Coxeter group with canonical generating set; in this paper we give an explicit description of the canonical generating set and the corresponding fundamental relations of $W^{\perp I}$ by using some graphs and some $2$-cell complex.
The conjugation action of the third factor $B_I$ on $W^{\perp I}$ induces automorphisms on $W^{\perp I}$ as a Coxeter system, i.e., those leaving the canonical generating set invariant.
We describe the subgroup $B_I$ as an extension of the fundamental group $Y_I$ of the above-mentioned $2$-cell complex by certain symmetries, and determine the action of $B_I$ in terms of this complex.

One may feel that the above-mentioned result is similar to the preceding result on the normalizer $N_W(W_I)$ of $W_I$ given by Brigitte Brink and Robert B.\ Howlett \cite{Bri-How} (or by Howlett \cite{How} for the case of a finite $W$).
In fact, our decomposition of the centralizer looks like the decomposition of the normalizer given by their work, and some properties of our decomposition are derived from the corresponding properties for normalizers.
However, it is worthy to emphasize that some fundamental properties of our decomposition are essentially new and desirable, and cannot be derived from the result on normalizers.
Moreover, those new properties enable us to prove (under a certain assumption on $I$) the following strong property, which will be included in a forthcoming paper by the author, of the centralizer $Z_W(W_I)$ that has been known for neither the centralizers nor the normalizers: The above-mentioned subgroup $Y_I$ of $Z_W(W_I)$ (which occupies a fairly large part of the third factor $B_I$) acts \emph{trivially} on every finite irreducible component of the Coxeter group $W^{\perp I}$.
(This novel property plays a central role in the author's recent study on the isomorphism problem in Coxeter groups; see \cite{Nui_ref}.)
Thus the result of this paper is not just an analogy of the preceding result, but gives further novel insights into the structural properties of Coxeter groups.
See Section \ref{sec:intro_related} for detailed explanations.
Conversely, as a by-product of our result, we also give another decomposition of $N_W(W_I)$ which shares the above-mentioned desirable properties with the decomposition of $Z_W(W_I)$.

\subsection{Related works}
\label{sec:intro_related}

As mentioned in Section \ref{sec:intro_background}, a similar problem of describing the normalizer $N_W(W_I)$ of a parabolic subgroup $W_I$ was studied for finite $W$ by Howlett \cite{How} and for general $W$ by Brink and Howlett \cite{Bri-How}.
Some part of our result is closely related to their preceding result, however some part of our result is essentially new and not derived from their work.
More precisely, in our description of the centralizer $Z_W(W_I)$ we use a groupoid $C$, which turns out to be a covering of the groupoid $G$ developed in their work, therefore some (but not all) properties of $C$ are inherited from those of $G$.
Our groupoid $C$ has $W^{\perp I} \rtimes Y_I$ as a vertex group, and by analyzing the groupoid we give a presentation of $W^{\perp I}$ as a Coxeter group and a description of $Y_I$ as the fundamental group of a $2$-cell complex $\mathcal{Y}$.
On the other hand, $N_W(W_I)$ is the semidirect product of $W_I$ and a vertex group of $G$, the latter being described in terms of a certain complex corresponding to $\mathcal{Y}$ (temporarily denoted by $\mathcal{Y}'$).
Now the preceding result gives a semidirect product decomposition of the vertex group of $G$, whose normal part is also a Coxeter group.
However, unlike the case of $Z_W(W_I)$ where the normal part $W^{\perp I}$ is a reflection subgroup of $W$, the normal part for the case of $N_W(W_I)$ is in general \emph{not} a reflection subgroup of $W$.
The property that our normal part $W^{\perp I}$ is a reflection subgroup plays a significant role in the proof of the above-mentioned result on the action of $Y_I$ on finite irreducible components of $W^{\perp I}$.
The properties of the decomposition $W^{\perp I} \rtimes Y_I$ are not immediately inherited from those of the vertex group of $G$ (see Remark \ref{rem:relationofCandG} for details).
Moreover, the structure of our complex $\mathcal{Y}$ is somewhat simpler than their complex $\mathcal{Y}'$; namely, $\mathcal{Y}$ has no $2$-cells with non-simple boundaries, while $\mathcal{Y}'$ may have such a non-simple cell.
Due to these differences, our description of the structure of $Z_W(W_I)$ is still hard to derive even if that of $N_W(W_I)$ is clear.

A description of normalizers of finite parabolic subgroups $W_I$ in Coxeter groups $W$ (in fact, also in extensions of $W$ by Coxeter graph automorphisms) was also studied by Richard E.\ Borcherds in \cite{Bor}.
Although that paper seems to focus on the normalizers as its main object, the result includes the case of some other subgroups of the normalizers such as the centralizers.
Borcherds's argument also has a flavor similar to that in the present paper; for example, his result also gives a semidirect product decomposition into a Coxeter group and some outer factor, the latter being expressed by using its classifying category.
However, the argument based on classifying categories resulted in geometric objects of higher dimensions, while the present paper requires only two-dimensional cell complexes to describe the full structure of centralizers.
This simplicity enables us to obtain more explicit expressions of centralizers and, as a consequence, to reveal further properties of the centralizers (see a forthcoming paper by the author) which do not immediately follow from the results of the present and Borcherds's papers.

On the other hand, for the centralizer $Z_W(W_I)$, a special case where $I$ consists of a single generator was studied and a similar decomposition was given by Brink \cite{Bri}, although explicit generators of the Coxeter group part ($W^{\perp I}$ in our notation) were not given.
Patrick Bahls and Mike Mihalik gave a complete description of $Z_W(W_I)$ if $W$ is an \lq\lq even'' Coxeter group, namely if the product of any two distinct generators has either an even or infinite order, as is found in their preprint and a recent book by Bahls \cite{Bah}.
Their approach is highly different from ours.
The commensurators of $W_I$ in $W$ and the centralizers of parabolic subgroups $A_X$ in Artin-Tits groups $A$ of certain types were examined by Luis Paris \cite{Par_commens,Par_centra}.
Eddy Godelle \cite{God_thesis,God_paper} described the normalizers, centralizers and related objects of $A_X$ in $A$, and revealed some new relationship between Artin-Tits groups and Coxeter groups.

The results of this paper are applied to the author's recent studies on the isomorphism problem in Coxeter groups \cite{Nui_ref} and some relevant topics in structural properties of Coxeter groups \cite{Nui_ext,Nui_indec}.
Moreover, recently Ivonne J.\ Ortiz \cite{Ort_err} used the result of (a preliminary version of) this paper for explicitly computing the lower algebraic $K$-theory of $\Gamma_3 = O^+(3,1) \cap GL(4,\mathbb{Z})$ (see also \cite{Ort_org}).

\subsection{Organization of this paper}
\label{sec:intro_organization}

The paper is organized as follows.
Section \ref{sec:preliminaries} is devoted to some preliminaries, and Section \ref{sec:groupoidC} summarizes the properties of the groupoid $C$ that are obtained by lifting the corresponding properties of the groupoid $G$.
Section \ref{sec:decompofC} investigates the groupoid $C$ further, enhancing its similarity with Coxeter groups noticed by Brink and Howlett in \cite{Bri-How}, and derive the decomposition of the vertex group $C_I=W^{\perp I} \rtimes Y_I$ as well as the presentations of $W^{\perp I}$ and $Y_I$.
(Although this looks similar to the decomposition $G_I=\widetilde{N}_I \rtimes M_I$ given in \cite{Bri-How}, there seems no overall relation between the structures of $\widetilde{N}_I$ and $W^{\perp I}$ or between those of $M_I$ and $Y_I$; see Remark \ref{rem:relationofCandG} for details.)
Section \ref{sec:descriptionoffullZ} completes the description of the entire $Z_W(W_I)$, and gives a description of $N_W(W_I)$ based on the same argument.
Section \ref{sec:example} carries out an explicit computation of the presentation of $Z_W(W_I)$ with an example.

%%%%
\paragraph*{Acknowledgments.}
This paper is an enhancement of the material in Master's thesis \cite{Nui} by the author.
The author would like to express his deep gratitude to everyone who helped him, especially to Professor Itaru Terada who was the supervisor of the author during the graduate course, and to Professor Kazuhiko Koike, for their invaluable advice and encouragement.
The author would also like to thank very much the anonymous referee for precious comments and suggestions to improve this paper, including separation of the final section of the original version as an individual paper. 
The present work was inspired by Brink's paper \cite{Bri} in the first place, and the paper \cite{Bri-How} of Brink and Howlett incited various improvements.
A part of this work was supported by JSPS Research Fellowship (No.\ 16-10825).
%%%%
%%%%%%%
\section{Preliminaries}
\label{sec:preliminaries}
%%%%%
\subsection{Groupoids and related objects}
\label{sec:groupoids}
We briefly summarize the notion of groupoids, fixing some notation and convention.
We refer to \cite{Coh} for the terminology and facts not mentioned here.

A \emph{groupoid} is a small category whose morphisms are all invertible; namely a family of sets $G=\{G_{x,y}\}_{x,y \in I}$ with index set $\mathcal{V}(G)=I$ endowed with (1) multiplications $G_{x,y} \times G_{y,z} \to G_{x,z}$ satisfying the associativity law, (2) an identity element $1=1_x$ in each $G_{x,x}$, and (3) an inverse $g^{-1} \in G_{y,x}$ for every $g \in G_{x,y}$.
In particular, each $G_{x,x}$ is a group, called a \emph{vertex group} of $G$.
We write $g \in G$ to signify $g \in \bigsqcup_{x,y \in I}G_{x,y}$, and $X \subseteq G$ for $X \subseteq \bigsqcup_{x,y \in I}G_{x,y}$.
If $g \in G_{y,x}$ we call $y$ and $x$ the \emph{target} and \emph{source} of $g$, respectively.
(Note that this is reverse to the usual convention.
The reason is that we will consider later an action of an element $g \in G_{y,x}$ of a certain groupoid $G$ from the left that maps an element associated to $x$ to one associated to $y$.)
An element $g \in G$ is called a \emph{loop} in $G$ if the source and target of $g$ are equal, namely if $g$ lies in some vertex group of $G$.
A sequence $(g_l,g_{l-1},\dots,g_1)$ of elements of $G$ is \emph{composable} if the source of $g_{i+1}$ is equal to the target of $g_i$ for every $1 \leq i \leq l-1$.
A \emph{homomorphism} between groupoids is a covariant functor between them regarded as categories.
Notions such as \emph{isomorphisms} and \emph{subgroupoids} are defined as usual.

Let $X \subseteq G$.
A \emph{word} in the set $X^{\pm 1}=X \sqcup \{x^{-1} \mid x \in X\}$ is a composable sequence $(x_l^{\varepsilon_l},\dots,x_2^{\varepsilon_2},x_1^{\varepsilon_1})$ of symbols in $X^{\pm 1}$.
We say that this word \emph{represents} the element $x_l^{\varepsilon_l} \cdots x_2^{\varepsilon_2}x_1^{\varepsilon_1}$ of $G$.
The set $X$ is called a \emph{generating set} of $G$ if every element of $G$ is represented by a word in $X^{\pm 1}$.
Given a generating set $X$ of $G$, a set $R$ of formal expressions of the form $\xi=\eta$, where $\xi$ and $\eta$ are words in $X^{\pm 1}$, is called a set of \emph{fundamental relations} of $G$ with respect to $X$ if (1) the two words $\xi$ and $\eta$ in each expression $\xi=\eta$ in $R$ represent the same element of $G$, and (2) if two words $\zeta$ and $\omega$ in $X^{\pm 1}$ represent the same element of $G$, then $\omega$ is obtained from $\zeta$ by a finite steps of transformations replacing a subword matching one side of an expression in $R \cup \{x^{\varepsilon}x^{-\varepsilon}=\emptyset \mid x \in X \textrm{ and } \varepsilon \in \{\pm 1\}\}$ (where $\emptyset$ represents the empty word) by the other side (or vice versa).

Let $\mathcal{G}$ be an unoriented graph on vertex set $\mathcal{V}$ with no loops.
The \emph{fundamental groupoid} of $\mathcal{G}$, denoted here by $\pi_1(\mathcal{G};*,*)$, is a groupoid with index set $\mathcal{V}$.
The edges of $\mathcal{G}$, each endowed with one of the two possible orientations, say, from $x \in \mathcal{V}$ to $y \in \mathcal{V}$, define elements of $\pi_1(\mathcal{G};y,x) = \pi_1(\mathcal{G};*,*)_{y,x}$, and $\pi_1(\mathcal{G};*,*)$ is freely generated by these elements.
Note that the same edge endowed with two opposite orientations give inverse elements to each other.
A vertex group $\pi_1(\mathcal{G};x)=\pi_1(\mathcal{G};x,x)$ is the \emph{fundamental group} of $\mathcal{G}$ at $x$, which is a free group.
To keep consistency with our convention for groupoids, in this paper \emph{we write a path in a graph from right to left}, hence a path $e_l \cdots e_2e_1$, where $e_i$ is an oriented edge from $x_{i-1}$ to $x_i$ for $1 \leq i \leq l$, determines an element of $\pi_1(\mathcal{G};x_l,x_0)$.

In this paper, the term \lq\lq \emph{complex}'' refers to a pair $\mathcal{K}=(\mathcal{K}^1,\mathcal{C})$ of a graph $\mathcal{K}^1$, called the \emph{$1$-skeleton} of $\mathcal{K}$, and a set $\mathcal{C}$ of some (not necessarily simple) closed paths in $\mathcal{K}$.
As an analogy of usual $2$-dimensional cell complexes, each closed path $c \in \mathcal{C}$ is regarded as a \lq\lq $2$-cell'' of $\mathcal{K}$ with boundary $c$.
When $\mathcal{K}^1$ has no loops, the \emph{fundamental groupoid} of $\mathcal{K}$, denoted similarly by $\pi_1(\mathcal{K};*,*)$, is constructed from $\pi_1(\mathcal{K}^1;*,*)$ by adding new relations \lq\lq boundary of a $2$-cell $=1$''.
Note that to every complex $\mathcal{K}$ a topological space (called a \emph{geometric realization} of $\mathcal{K}$) is associated in such a way that $\pi_1(\mathcal{K};*,*)$ is naturally a full subgroupoid of the fundamental groupoid of the space (see \cite[Section 5.6]{Coh} for details).

A groupoid $\widetilde{G}$ is called a \emph{covering groupoid} of a groupoid $G$ with \emph{covering map} $p:\widetilde{G} \to G$ if $p$ is a surjective homomorphism and, moreover, for any $g \in G$ with source $x$ and $\widetilde{x} \in \mathcal{V}(\widetilde{G})$ projecting onto $x$ (called a \emph{lift} of $x$), there exists a unique element $\widetilde{g} \in \widetilde{G}$ with source $\widetilde{x}$ projecting onto $g$ (also called a \emph{lift} of $g$); see e.g., \cite[Section 9.2]{Bro}.
Any lift of an identity element in $G$ is an identity element in $\widetilde{G}$.
If $(g_l,\dots,g_2,g_1)$ is a composable sequence of elements of $G$, then for any lift $\widetilde{x}$ of the source $x$ of $g_1$, there is a unique composable sequence $(\widetilde{g_l},\dots,\widetilde{g_2},\widetilde{g_1})$ of lifts of $g_l,\dots,g_2,g_1$ respectively, with $\widetilde{x}$ being the source of $\widetilde{g_1}$; and the product $\widetilde{g_l} \cdots \widetilde{g_2}\widetilde{g_1}$ is the lift of $g_l \cdots g_2g_1$ with source $\widetilde{x}$.
In later sections we will use the following property of covering groupoids, whose proof is straightforward and so omitted here:
\begin{prop}
\label{prop:liftofrelations}
Suppose that a groupoid $G$ admits a covering groupoid $\widetilde{G}$, a generating set $X$ and a set $R$ of fundamental relations.
Then the set $\widetilde{X}$ of the lifts of the generators $g \in X$ generates $\widetilde{G}$.
Moreover, $\widetilde{G}$ admits, with respect to the $\widetilde{X}$, a set $\widetilde{R}$ of fundamental relations consisting of the formal expressions $\widetilde{\xi}=\widetilde{\eta}$ such that $\widetilde{\xi}$ and $\widetilde{\eta}$ are words in ${\widetilde{X}}^{\pm 1}$ with the same source and these are lifts of corresponding terms of an expression $\xi=\eta$ belonging to $R$.
\end{prop}
%%
%%%%%
\subsection{Coxeter groups}
\label{sec:Coxetergroups}
The basics of Coxeter groups summarized here are found in \cite{Hum} unless otherwise noticed.
%%%
\subsubsection{Definitions}
\label{sec:defofCox}
A pair $(W,S)$ of a group $W$ and its generating set $S$ is called a \emph{Coxeter system} if $W$ admits the following presentation
\begin{displaymath}
W=\langle S \mid (st)^{m_{s,t}}=1 \mbox{ for all } s,t \in S \mbox{ such that } m_{s,t}<\infty \rangle \enspace,
\end{displaymath}
where the $m_{s,t} \in \{1,2,\dots\} \cup \{\infty\}$ are symmetric in $s$ and $t$, and $m_{s,t}=1$ if and only if $s=t$.
Note that \emph{the set $S$ may have infinite cardinality}.
A group $W$ is called a \emph{Coxeter group} if $(W,S)$ is a Coxeter system for some $S \subseteq W$.
An \emph{isomorphism of Coxeter systems} from $(W,S)$ to $(W',S')$ is a group isomorphism $W \to W'$ that maps $S$ onto $S'$.
It is shown that $m_{s,t}$ is precisely the order of $st \in W$, therefore the system $(W,S)$ determines uniquely the data $(m_{s,t})_{s,t}$ and hence the \emph{Coxeter graph} $\Gamma$, which is a simple unoriented graph with vertex set $S$ in which two vertices $s,t \in S$ are joined by an edge with label $m_{s,t}$ if and only if $m_{s,t} \geq 3$ (by usual convention, the label is omitted when $m_{s,t}=3$).
For $I \subseteq S$, the subgroup $W_I = \langle I \rangle$ of $W$ generated by $I$ is called a \emph{parabolic subgroup}.
If $I$ is (the vertex set of) a connected component of $\Gamma$, then $W_I$ and $I$ are called an \emph{irreducible component} of $W$ (or of $(W,S)$, if we emphasize the generating set $S$) and of $S$, respectively.
Now $W$ is the (restricted) direct product of its irreducible components.
If $\Gamma$ is connected, then $(W,S)$, $W$ and $S$ are called \emph{irreducible}.
If $I \subseteq S$, then $(W_I,I)$ is also a Coxeter system, of which the Coxeter graph $\Gamma_I$ is the full subgraph of $\Gamma$ with vertex set $I$ and the length function is the restriction of the length function $\ell$ of $(W,S)$.
%%%%%
\subsubsection{Geometric representation and root systems}
\label{sec:rootsystem}
Let $V$ denote the space of the \emph{geometric representation} of $(W,S)$ over $\mathbb{R}$ equipped with a basis $\Pi=\{\alpha_s \mid s \in S\}$ and a $W$-invariant symmetric bilinear form $\langle \,,\, \rangle$ determined by
\begin{displaymath}
\langle \alpha_s, \alpha_t \rangle =
\begin{cases}
-\cos(\pi / m_{s,t}) & \mbox{if } m_{s,t} < \infty \enspace; \\
-1 & \mbox{if } m_{s,t} = \infty \enspace,
\end{cases}
\end{displaymath}
where $W$ acts faithfully on $V$ by $s \cdot v=v-2\langle \alpha_s, v\rangle \alpha_s$ for $s \in S$ and $v \in V$.
Then the \emph{root system} $\Phi=W \cdot \Pi$ consists of unit vectors with respect to the $\langle \,,\, \rangle$, and it is the disjoint union $\Phi=\Phi^+ \sqcup \Phi^-$ of $\Phi^+=\Phi \cap \mathbb{R}_{\geq 0}\Pi$ and $\Phi^-=-\Phi^+$ where $\mathbb{R}_{\geq 0}\Pi$ signifies the set of nonnegative linear combinations of elements of $\Pi$.
Elements of $\Phi$, $\Phi^+$, and $\Phi^-$ are called \emph{roots}, \emph{positive roots}, and \emph{negative roots}, respectively.
For any subset $\Psi \subseteq \Phi$ and $w \in W$, put
\begin{displaymath}
\Psi^{\pm}=\Psi \cap \Phi^{\pm}\,, \mbox{ respectively,}
\end{displaymath}
and
\begin{displaymath}
\Psi[w] =\{\gamma \in \Psi^+ \mid w \cdot \gamma \in \Phi^-\} \enspace.
\end{displaymath}
We have $\ell(w)=|\Phi[w]|$, hence $w=1$ if and only if $\Phi[w]=\emptyset$.
This implies that the set $\Phi[w]$ determines the element $w \in W$ uniquely.
The following property is well known:
\begin{lem}
\label{lem:lengthofmultiple}
For $w_1,w_2 \in W$, the following conditions are equivalent:
\begin{enumerate}
\item $\ell(w_1w_2)=\ell(w_1)+\ell(w_2)$;
\item $\Phi[w_2] \subseteq \Phi[w_1w_2]$;
\item $\Phi[w_1] \cap \Phi[w_2^{-1}]=\emptyset$;
\item $\Phi[w_1w_2]=(w_2^{-1} \cdot \Phi[w_1]) \sqcup \Phi[w_2]$.
\end{enumerate}
\end{lem}

For any element $v=\sum_{s \in S}c_s\alpha_s \in V$, define the \emph{support} of $v$ by
\begin{displaymath}
\mathrm{Supp}\,v=\{s \in S \mid c_s \neq 0\} \enspace.
\end{displaymath}
For each $I \subseteq S$, put
\begin{displaymath}
\Pi_I=\{\alpha_s \mid s \in I\} \subseteq \Pi \,,\, V_I=\mathrm{Supp}\,\Pi_I \subseteq V \mbox{ and } \Phi_I=\Phi \cap V_I \enspace.
\end{displaymath}
It is well known that $\Phi_I$ coincides with the root system $W_I \cdot \Pi_I$ of $(W_I,I)$ (see also Theorem \ref{thm:reflectionsubgroup_Deodhar}(3) below).
The support of any $\gamma \in \Phi$ is irreducible; this follows from the facts that we have $\gamma \in \Phi_I$ where $I= \mathrm{Supp}\,\gamma$, and that $\mathrm{Supp}\, (w \cdot \alpha_s)$ (for $w \in W_I$ and $s \in I$) is contained in the irreducible component of $I$ containing $s$.
%%%%%
\subsubsection{Reflection subgroups}
\label{sec:reflectionsubgroup}
For a $\gamma=w \cdot \alpha_s \in \Phi$, let $s_\gamma=wsw^{-1}$ denote the \emph{reflection} along $\gamma$ acting on $V$ by $s_\gamma \cdot v=v-2 \langle \gamma, v \rangle \gamma$ for $v \in V$.
Note that $\ell(us_\gamma) \neq \ell(u)$ for any $u \in W$.
For any subset $\Psi \subseteq \Phi$, let $W(\Psi)$ denote the subgroup generated by all $s_\gamma$ with $\gamma \in \Psi$.
A subgroup of $W$ of this form is called a \emph{reflection subgroup}.
Note that $W(\Psi)=W(W(\Psi) \cdot \Psi)$ and $-W(\Psi) \cdot \Psi = W(\Psi) \cdot \Psi$.
Moreover, let $\Pi(\Psi)$ denote the set of all \lq\lq simple roots'' $\gamma \in (W(\Psi) \cdot \Psi)^+$, namely all the $\gamma$ such that any expression $\gamma=\sum_{i=1}^{r}c_i\beta_i$ with $c_i>0$ and $\beta_i \in (W(\Psi) \cdot \Psi)^+$ satisfies that $\beta_i=\gamma$ for all $i$.
Put
\begin{displaymath}
S(\Psi)=\{s_\gamma \mid \gamma \in \Pi(\Psi)\} \enspace.
\end{displaymath}
It was shown by Vinay V.\ Deodhar \cite{Deo_refsub} that $(W(\Psi),S(\Psi))$ is a Coxeter system.
Note that Matthew Dyer \cite{Dye} also proved independently that $W(\Psi)$ is a Coxeter group and determined a corresponding generating set, which in fact coincides with $S(\Psi)$.
First, we summarize Deodhar's result (the claim 1 below) and some related properties:
\begin{thm}
\label{thm:reflectionsubgroup_Deodhar}
Let $\Psi \subseteq \Phi$ be any subset.
\begin{enumerate}
\item The pair $(W(\Psi),S(\Psi))$ is a Coxeter system.
For any $\beta \in \Pi(\Psi)$, the reflection $s_\beta$ permutes the elements of $(W(\Psi) \cdot \Psi)^+ \smallsetminus \{\beta\}$.
Moreover, any $\gamma \in (W(\Psi) \cdot \Psi)^+$ can be expressed as a nonnegative linear combination of elements of $\Pi(\Psi)$.
\item Let $\beta,\gamma \in W(\Psi) \cdot \Psi$ and $w \in W(\Psi)$.
Then $ws_\beta w^{-1}=s_\gamma$ if and only if $w \cdot \beta=\pm \gamma$.
\item If $\gamma \in \Phi$ and $s_\gamma \in W(\Psi)$, then $\gamma \in W(\Psi) \cdot \Pi(\Psi)$.
In particular, we have $W(\Psi) \cdot \Psi = W(\Psi) \cdot \Pi(\Psi)$.
\item Let $\ell_\Psi$ be the length function of $(W(\Psi),S(\Psi))$.
Then for $w \in W(\Psi)$ and $\gamma \in (W(\Psi) \cdot \Psi)^+$, we have $\ell_\Psi(ws_\gamma)<\ell_\Psi(w)$ if and only if $w \cdot \gamma \in \Phi^-$.
\item If $w \in W(\Psi)$, then $\ell_\Psi(w)=|(W(\Psi) \cdot \Psi)[w]|<\infty$.
\end{enumerate}
\end{thm}
Note that the claim 2 is well known for arbitrary roots.
On the other hand, once the claim 3 is proven, the claims 4 and 5 will follow from the proof of the main theorem of Deodhar's another paper \cite{Deo_char}.
To prove the claim 3, we require another theorem of Deodhar \cite{Deo_root}:
\begin{thm}
[{\cite[Theorem 5.4]{Deo_root}}]
\label{thm:involution}
If $w \in W$ and $w^2=1$, then $w$ decomposes into a product of commuting reflections along pairwise orthogonal roots.
Hence if $w$ is an involution, then $w \cdot \gamma=-\gamma$ for some $\gamma \in \Phi$.
\end{thm}
\begin{proof}
[Proof of Theorem \ref{thm:reflectionsubgroup_Deodhar}(3)]
For the first part, let $\gamma \in \Phi$ such that $s_\gamma \in W(\Psi)$.
Then Theorem \ref{thm:involution} gives us a decomposition $s_\gamma=s_{\beta_1} \dotsm s_{\beta_r}$ of the involution $s_\gamma$ in a Coxeter group $W(\Psi)$ into pairwise commuting distinct reflections $s_{\beta_i}$ along $\beta_i \in W(\Psi) \cdot \Pi(\Psi)$.
Write $\gamma=w \cdot \alpha_s$ with $s \in S$ and $w \in W$, then $s=w^{-1}s_\gamma w=s_{\beta'_1} \dotsm s_{\beta'_r}$ where all $\beta'_i=w^{-1} \cdot \beta_i$ are orthogonal as well as the $\beta_i$.
This implies that $s \cdot \beta'_i=-\beta'_i$, therefore $\beta'_i=\pm \alpha_s$ for every $i$.
Hence we have $r=1$ and $\gamma=\pm w \cdot \beta'_1=\pm \beta_1 \in W(\Psi) \cdot \Pi(\Psi)$, as desired.

For the second part, any element $\gamma \in W(\Psi) \cdot \Psi$ satisfies that $s_{\gamma} \in W(W(\Psi) \cdot \Psi) = W(\Psi)$, therefore $\gamma \in W(\Psi) \cdot \Pi(\Psi)$ by the first part.
The other inclusion is obvious.
\end{proof}
Secondly, we summarize a part of Dyer's result and its consequences needed later.
We say that a subset $\Psi \subseteq \Phi^+$ is a \emph{root basis} if for all $\beta,\gamma \in \Psi$, we have
\begin{displaymath}
\begin{cases}
\langle \beta,\gamma \rangle=-\cos(\pi/m) & \mbox{if } s_\beta s_\gamma \mbox{ has order } m<\infty \enspace;\\
\langle \beta,\gamma \rangle \leq -1 & \mbox{if } s_\beta s_\gamma \mbox{ has infinite order}.
\end{cases}
\end{displaymath}
\begin{thm}
[{\cite[Theorem 4.4]{Dye}}]
\label{thm:conditionforrootbasis}
Let $\Psi \subseteq \Phi^+$ be any subset.
Then we have $\Pi(\Psi)=\Psi$ if and only if $\Psi$ is a root basis.
\end{thm}
\begin{cor}
\label{cor:fintyperootbasis}
Let $\Psi \subseteq \Phi^+$ be a root basis such that $|W(\Psi)|<\infty$.
Then $\Psi$ is a basis of a positive definite subspace of $V$ with respect to the form $\langle \,,\, \rangle$.
\end{cor}
\begin{proof}
By Theorem \ref{thm:conditionforrootbasis}, the finite Coxeter group $W(\Psi)$ has its own geometric representation space $V'$ with basis $\Pi'=\{\overline{\beta} \mid \beta \in \Psi\}$ and positive definite bilinear form $\langle \,,\, \rangle'$ (see \cite[Theorem 6.4]{Hum}) such that $\langle \beta, \gamma \rangle=\langle \overline{\beta}, \overline{\gamma} \rangle'$ for $\beta,\gamma \in \Psi$.
Now if $v=\sum_{\beta \in \Psi}c_\beta \beta \in \mathrm{span}\,\Psi$ and $\langle v, v \rangle \leq 0$, then $\langle \overline{v}, \overline{v} \rangle'=\langle v, v \rangle \leq 0$ where $\overline{v}=\sum_{\beta \in \Psi}c_\beta \overline{\beta} \in V'$, implying that $c_\beta=0$ as desired.
\end{proof}
Here we present two more results needed later as well.
The first one is a slight improvement of \cite[Proposition 2.6]{How-Row-Tay}.
Note that the assumption $|S|<\infty$ in the original statement is in fact not necessary; namely we have:
\begin{prop}
\label{prop:finitesubsystem}
Let $\Psi \subseteq \Phi^+$ be a root basis such that $|W(\Psi)|<\infty$, and $U= \mathrm{span}\,\Psi$.
Then there exist $w \in W$ and $I \subseteq S$ such that $|W_I|<\infty$ and $w \cdot (U \cap \Phi^+)=\Phi_I^+$.
Moreover, this $w$ maps $U \cap \Pi$ into $\Pi_I$.
\end{prop}
\begin{proof}
First, Theorem \ref{thm:reflectionsubgroup_Deodhar}(3) implies that the set $\{\gamma \in \Phi \mid s_\gamma \in W(\Psi)\}$ coincides with the root system $W(\Psi) \cdot \Psi$ of $W(\Psi)$, which is finite by the hypothesis.
Thus the hypothesis of the original proposition \cite[Proposition 2.6]{How-Row-Tay} is satisfied, hence it follows from that proposition that there are $w \in W$ and $I \subseteq S$ such that $|W_I|<\infty$ and $w \cdot (U \cap \Phi)=\Phi_I$ (note that $U= \mathrm{span}\,(W(\Psi) \cdot \Psi)$).
Choosing such a $w$ with shortest length, it also holds that $w^{-1} \cdot \Phi_I^+ \subseteq U \cap \Phi^+$.
Indeed, if $\gamma \in \Phi_I[w^{-1}]$, then $s_\gamma w \in W$ also satisfies the condition and is shorter than $w$ (apply Theorem \ref{thm:reflectionsubgroup_Deodhar}(4) to $w^{-1}$).
Hence the first claim holds.

For the second claim, if $\alpha_s \in U \cap \Pi$ and $w \cdot \alpha_s=\sum_{t \in I}c_t\alpha_t$ with $c_t \geq 0$, then $\alpha_s=\sum_{t \in I}c_tw^{-1} \cdot \alpha_t$, while $w^{-1} \cdot \alpha_t \in \Phi^+$.
Thus we have $w^{-1} \cdot \alpha_t=\alpha_s$ whenever $c_t>0$, showing that $w \cdot \alpha_s=\alpha_t \in \Pi_I$ as desired.
\end{proof}
The second result is that the theorem \cite[Theorem 1.12 (d)]{Hum} also holds for a general $W$:
\begin{prop}
\label{prop:fixingroots}
Let $\Psi \subseteq \Phi$ be any subset.
Then any $w \in W$ of finite order that fixes $\Psi$ pointwise decomposes into a product of reflections also having this property.
\end{prop}
\begin{proof}
First, a well-known theorem of Jacques Tits implies that the finite subgroup $\langle w \rangle$ of $W$ is conjugate to a subgroup of a finite $W_I$ (see \cite[Theorem 4.5.3]{Bjo-Bre} for a proof), therefore we may assume without loss of generality that $w \in W_I$.
Since $V_I$ is positive definite (see \cite[Theorem 6.4]{Hum}), each $\gamma \in \Psi$ decomposes as $\gamma=v+v'$ with $v \in V_I$ and $v' \in {V_I}^\perp$.
Now $w \in W_I$ fixes every $v=\gamma-v'$, therefore the result for the finite $W_I$ yields a decomposition $w=s_{\beta_1} \cdots s_{\beta_r}$, where each reflection $s_{\beta_j} \in W_I$ fixes every $v$, hence fixes every $\gamma=v+v' \in \Psi$, as desired.
\end{proof}
%%
%%%%%%
\subsubsection{Finite parabolic subgroups and their longest elements}
\label{sec:longestelement}
We say that a subset $I \subseteq S$ is of \emph{finite type} if $|W_I|<\infty$, or equivalently $|\Phi_I|<\infty$.
The finite irreducible Coxeter groups have been classified as summarized in \cite[Chapter 2]{Hum}.
Here we give a canonical labelling $r_1,r_2,\dots,r_n$ (where $n = |I|$) of elements of an irreducible subset $I \subseteq S$ of each finite type used in this paper in the following manner, where we put $m_{i,j}=m_{r_i,r_j}$ for simplicity and the $m_{i,j}$ not listed here are all equal to $2$:
\begin{description}
\item[Type $A_n$ ($n \geq 2$):] $m_{i,i+1}=3$ ($1 \leq i \leq n-1$);
\item[Type $B_n$ ($n \geq 2$):] $m_{i,i+1}=3$ ($1 \leq i \leq n-2$) and $m_{n-1,n}=4$;
\item[Type $D_n$ ($n \geq 4$):] $m_{i,i+1}=m_{n-2,n}=3$ ($1 \leq i \leq n-2$);
\item[Type $E_n$ ($n=6,7,8$):] $m_{1,3}=m_{2,4}=m_{i,i+1}=3$ ($3 \leq i \leq n-1$);
\item[Type $F_4$:] $m_{1,2}=m_{3,4}=3$ and $m_{2,3}=4$;
\item[Type $H_n$ ($n=3,4$):] $m_{1,2}=5$ and $m_{i,i+1}=3$ ($2 \leq i \leq n-1$);
\item[Type $I_2(m)$ ($m \geq 5$):] $m_{1,2}=m$.
\end{description}

Let $w_0(I)$ denote the (unique) longest element of a finite $W_I$, which has order two and maps $\Pi_I$ onto $-\Pi_I$.
Now let $I$ be irreducible of finite type.
If $I$ is of type $A_n$ ($n \geq 2$), $D_k$ ($k$ odd), $E_6$ or $I_2(m)$ ($m$ odd), then the automorphism of the Coxeter graph $\Gamma_I$ of $W_I$ induced by (the conjugation action of) $w_0(I)$ is the unique nontrivial automorphism on $\Gamma_I$.
Otherwise $w_0(I)$ lies in the center $Z(W_I)$ of $W_I$ and the induced automorphism on $\Gamma_I$ is trivial, in which case we say that $I$ is of \emph{$(-1)$-type}.
Moreover, if $W_I$ is finite but not irreducible, then $w_0(I)=w_0(I_1) \dotsm w_0(I_k)$ where the $I_i$ are the irreducible components of $I$.
Note that for an arbitrary $I \subseteq S$, the center $Z(W_I)$ is an elementary abelian $2$-group generated by the $w_0(J)$ where $J$ runs over all irreducible components of $I$ of $(-1)$-type.
%%%%%
\section{The groupoid $C$}
\label{sec:groupoidC}

From now on, we fix an arbitrarily given subset $I \subseteq S$ unless specifically noted otherwise.
In this section, we introduce and study a groupoid $C$, one of whose vertex groups gives a fairly large part (still not all) of the centralizer $Z_W(W_I)$.
As we shall see in Section \ref{sec:groupoidC_lifting}, our groupoid $C$ is a covering of the groupoid $G$ defined by Brink and Howlett in \cite{Bri-How} for their study of the normalizer $N_W(W_I)$.
As a consequence, many properties, but not all, of $C$ can be obtained from those of $G$.

Here we survey the connection among the centralizer $Z_W(W_I)$ and the two groupoids $C$ and $G$.
If $X$ is a $W$-set (that is, a set with $W$-action) and $Y \subseteq X$, let $\mathcal{G}=(X;Y)$ denote tentatively a groupoid with index set $Y$ defined by $\mathcal{G}=\{\mathcal{G}_{y,y'}\}_{y,y' \in Y}$ and $\mathcal{G}_{y,y'}=\{w \in W \mid y=w \cdot y'\}$, with multiplication induced by that of $W$.
The vertex group $\mathcal{G}_{y,y}$ of $\mathcal{G}$ at $y \in Y$ is the stabilizer of $y$ in $W$.
If $X'$ is another $W$-set and $Y' \subseteq X'$, then a $W$-equivalent map $\varphi:X \to X'$ with $\varphi(Y) \subseteq Y'$ induces a groupoid homomorphism $\varphi_*:(X;Y) \to (X';Y')$.
It is a covering of groupoids if $\varphi(Y)=Y'$ and $Y=\varphi^{-1}(Y')$.

Now fix a set $\Lambda$ with $|\Lambda|=|I|$.
For a $W$-set $X$, let $X^*$ temporarily denote the set of all injections $\Lambda \to X$ with the $W$-action induced from that on $X$.
Then the centralizer $Z_W(W_I)$, that is the final target of our study, is the vertex group of a groupoid $(T^*;S^*)$ at an element $x_I \in S^*$, where $T$ denotes the set of reflections in $W$ with respect to $S$ (on which $W$ acts by conjugation) and the image of the map $x_I$ is $I$.

Look at the following diagram of groupoids, whose morphisms are induced by the obvious $W$-equivalent maps.
Here the symbol $\binom{X}{|\Lambda|}$, where $X$ is a $W$-set, denotes the set of all subsets of $X$ with cardinality $|\Lambda|$ with the $W$-action induced from that on $X$.
The two morphisms marked $\ast$ are coverings, and the one marked $\#$ is an inclusion as a full subgroupoid.
Note that the composition of the vertical arrows induces a bijection on the index sets.
The groupoid $C$ we introduce below is isomorphic to the connected component of the groupoid $(\Phi^*;\Pi^*)$ on the top-left corner containing the lift of $x_I$.
Since $Z_W(W_I)$ contains elements projected from those of the groupoid $(\Phi^*;(\pm \Pi)^*)$ with source being the lift of $x_I \in S^*$ that lies in $\Pi^* \subseteq (\pm \Pi)^*$ and target being its other lifts in $(\pm \Pi)^* \smallsetminus \Pi^*$, the vertex group $C_I$ is only a part of $Z_W(W_I)$ in general.
On the other hand, the groupoid $G$ in \cite{Bri-How} is a connected component of the groupoid $(\binom{\Phi}{|\Lambda|};\binom{\Pi}{|\Lambda|})$ on the right.
\begin{displaymath}
\begin{CD}
\hspace*{2.5em} C_I \quad \subseteq \quad (\Phi^*;\Pi^*) @>*>> (\binom{\Phi}{|\Lambda|};\binom{\Pi}{|\Lambda|}) \quad \supseteq \quad G\\
\hspace*{7em} @V\#VV\\
\hspace*{7em} (\Phi^*;(\pm \Pi)^*)\\
\hspace*{7em} @V*VV\\
Z_W(W_I) \quad \subseteq \quad (T^*;S^*)
\end{CD}
\end{displaymath}
%%%%%
\subsection{Definitions}
\label{sec:groupoidC_def}
In what follows, we would like to deal with an \lq\lq ordered tuple'' consisting of the elements of a given subset of $S$.
The reason is that elements of the centralizer $Z_W(W_I)$ not only leaves the set $I$ invariant but also fix every element of $I$ pointwise.
For the purpose, first we fix an auxiliary index set $\Lambda_0$ having the same cardinality as $S$.
Then put
\begin{displaymath}
S^{(\Lambda)} = \{x:\Lambda \to S \mid x \mbox{ is injective}\} \mbox{ for each } \Lambda \subseteq \Lambda_0
\end{displaymath}
and
\begin{displaymath}
S^* = \bigcup_{\Lambda \subseteq \Lambda_0} S^{(\Lambda)} \enspace.
\end{displaymath}
For $x \in S^{(\Lambda)}$ and $\lambda \in \Lambda$, we write $x_\lambda$ for $x(\lambda)$; thus $x$ may be regarded as a duplicate-free \lq\lq $\Lambda$-tuple'' $(x_\lambda)=(x_\lambda)_{\lambda \in \Lambda}$ of elements of $S$.
If $|\Lambda| = n < \infty$, then it is regarded as a duplicate-free sequence $(x_1,x_2,\dots,x_n)$ by identifying $\Lambda$ with $\{1,2,\dots,n\}$.
We write
\begin{displaymath}
x_A = \{x_\lambda \mid \lambda \in A\} \mbox{ for each $A \subseteq \Lambda$\,, and } [x] = x_{\Lambda} \enspace.
\end{displaymath}

If we put
\begin{displaymath}
\widehat{C}_{x,y}=\{w \in W \mid \alpha_{x_\lambda}=w \cdot \alpha_{y_\lambda} \mbox{ for all } \lambda \in \Lambda\} \mbox{ for all } x,y \in S^{(\Lambda)}
\end{displaymath}
and $\widehat{C}_{x,y} = \emptyset$ for all $x \in S^{(\Lambda)}$ and $y \in S^{(\Lambda')}$ such that $\Lambda \neq \Lambda'$, then $\widehat{C}=\{\widehat{C}_{x,y}\}_{x,y \in S^*}$ is a groupoid with vertex set $\mathcal{V}(\widehat{C})=S^*$ and multiplication induced by that of $W$.
Write
\begin{equation}
\label{eq:defofwcdoty}
x = w \cdot y \mbox{ if } x,y \in S^* \mbox{ and } w \in \widehat{C}_{x,y} \enspace.
\end{equation}
Now we fix elements
\begin{displaymath}
\Lambda \subseteq \Lambda_0 \mbox{ and } x_I \in S^{(\Lambda)} \mbox{ such that } [x_I] = I \enspace,
\end{displaymath}
and define the groupoid $C$ as the connected component of $\widehat{C}$ containing $x_I$; hence
\begin{displaymath}
\mathcal{V}(C)=\{x \in S^{(\Lambda)} \mid \widehat{C}_{x,x_I} \neq \emptyset\} \mbox{ and } C_{x,y}=\widehat{C}_{x,y} \mbox{ for } x,y \in \mathcal{V}(C) \enspace.
\end{displaymath}
Note that $C_{x,y} \neq \emptyset$ for all $x,y \in \mathcal{V}(C)$.
We put
\begin{displaymath}
C_I=C_{x_I,x_I}=\{w \in W \mid w \cdot \alpha_s=\alpha_s \mbox{ for all } s \in I\} \enspace,
\end{displaymath}
which is a normal subgroup of $Z_W(W_I)$.
We will see below that this $C_I$ occupies a fairly large part, but in general not the whole, of $Z_W(W_I)$.

We will see in Theorem \ref{thm:decompofC} the following semidirect product decomposition
\begin{equation}
\label{eq:decompofC_I}
C_I=W^{\perp I} \rtimes Y_I \enspace.
\end{equation}
We define the factors below.
First, for arbitrary subsets $I,J \subseteq S$, let
\begin{displaymath}
\Phi_J^{\perp I}=\{\gamma \in \Phi_J \mid \mbox{$\gamma$ is orthogonal to every $\alpha_s \in \Pi_I$}\}
\end{displaymath}
and
\begin{displaymath}
W_J^{\perp I}=W(\Phi_J^{\perp I})
\end{displaymath}
(see Section \ref{sec:reflectionsubgroup} for notations).
Note that for any $\gamma \in \Phi_J^{\perp I}$, we have $s_{\gamma} \cdot \alpha_s = \alpha_s$ for every $s \in I$, therefore $s_{\gamma} \cdot \Phi_J^{\perp I} = \Phi_J^{\perp I}$.
Hence by Theorem \ref{thm:reflectionsubgroup_Deodhar}, $W_J^{\perp I}$ is a Coxeter group with root system $W_J^{\perp I} \cdot \Phi_J^{\perp I} = \Phi_J^{\perp I}$, generating set and simple system given by
\begin{displaymath}
R^{J,I}=S(\Phi_J^{\perp I}) \mbox{ and } \Pi^{J,I}=\Pi(\Phi_J^{\perp I})\,, \mbox{ respectively.}
\end{displaymath}
In the notations, the symbol $J$ will be omitted when $J=S$.
In particular, the factor $W^{\perp I}$ in (\ref{eq:decompofC_I}) is given by
\begin{displaymath}
W^{\perp I}=W_S^{\perp I}=\langle \{s_\gamma \mid \gamma \in \Phi^{\perp I}\} \rangle \enspace.
\end{displaymath}
For the factor $Y_I$ in (\ref{eq:decompofC_I}), first define a subgroupoid $\widehat{Y}=\{\widehat{Y}_{x,y}\}_{x,y \in S^*}$ of $\widehat{C}$ by
\begin{displaymath}
\widehat{Y}_{x,y}=\{w \in \widehat{C}_{x,y} \mid w \cdot (\Phi^{\perp [y]})^+ \subseteq \Phi^+\}=\{w \in \widehat{C}_{x,y} \mid (\Phi^{\perp [x]})^+=w \cdot (\Phi^{\perp [y]})^+\}
\end{displaymath}
for $x,y \in S^*$.
Here the second equality is a consequence of the following fact
\begin{displaymath}
\Phi^{\perp [y]}=w \cdot \Phi^{\perp [x]} \mbox{ for all } w \in \widehat{C}_{x,y} \enspace,
\end{displaymath}
also implying that $W^{\perp [x]}$ is normal in $\widehat{C}_{x,x}$.
Now let $Y$ be the full subgroupoid of $\widehat{Y}$ with vertex set $\mathcal{V}(C)$, namely
\begin{displaymath}
\mathcal{V}(Y)=\mathcal{V}(C) \mbox{ and } Y_{x,y}=\widehat{Y}_{x,y} \mbox{ for } x,y \in \mathcal{V}(C) \enspace,
\end{displaymath}
and write
\begin{displaymath}
Y_I = Y_{x_I,x_I} = \{w \in C_I \mid (\Phi^{\perp I})^+ = w \cdot (\Phi^{\perp I})^+\} \enspace.
\end{displaymath}
Then $Y$ is a subgroupoid of $C$.
Note that $W^{\perp [x]} \cap Y_{x,x}=1$ since any element of $W^{\perp [x]} \cap Y_{x,x}$ has length zero by virtue of Theorem \ref{thm:reflectionsubgroup_Deodhar}(5).
%%%%%
\subsection{Preceding results on normalizers}
\label{sec:precedingonN}
Here we summarize some preceding results on the normalizers $N_W(W_I)$ of $W_I$ in $W$ given by Brink and Howlett \cite{Bri-How}.
Their first step is the decomposition $N_W(W_I)=W_I \rtimes N_I$ (see \cite[Proposition 2.1]{Bri-How}), where the factor
\begin{displaymath}
N_I=\{w \in W \mid w \cdot \Pi_I=\Pi_I\}
\end{displaymath}
is the main subject of the paper \cite{Bri-How}.

The following theorem, used in \cite{Bri-How} for general $S$, is proven in \cite{Deo_root} by Deodhar under the assumption $|S|<\infty$.
A proof for the general case is given in \cite{Nui_indec}.
\begin{thm}
[Deodhar]
\label{thm:Deodhar_infinitecase}
Let $K \subseteq S$ be any subset and $J \subset K$ a proper subset, and suppose that $W_K$ is infinite and irreducible.
Then $|\Phi_K \smallsetminus \Phi_J|=\infty$.
\end{thm}
For $J,K \subseteq S$, let $J_{\sim K}$ denote the union of the connected components of the graph $\Gamma_{J \cup K}$ having nonempty intersection with $K$.
Then Theorem \ref{thm:Deodhar_infinitecase} implies that, in the case $J \cap K=\emptyset$, the set $J_{\sim K} \subseteq S$ is of finite type if and only if the set $\Phi_{J \cup K} \smallsetminus \Phi_J=\Phi_{J_{\sim K}} \smallsetminus \Phi_J$ is finite.
Here the last equality follows from the irreducibility of the support of any root (see Section \ref{sec:rootsystem}).
In this case, it is shown in \cite{Bri-How} that
\begin{displaymath}
w_0(J_{\sim K})w_0(J_{\sim K} \smallsetminus K) \cdot \Pi_J \subseteq \Pi_{J \cup K}
\end{displaymath}
and
\begin{displaymath}
\Phi\left[w_0(J_{\sim K})w_0(J_{\sim K} \smallsetminus K)\right]=\Phi_{J \cup K}^+ \smallsetminus \Phi_J \enspace.
\end{displaymath}
Moreover, it is easy to show that
\begin{equation}
\label{eq:inverseofw0w0}
(w_0(J_{\sim K})w_0(J_{\sim K} \smallsetminus K))^{-1}=w_0(J'_{\sim K'})w_0(J'_{\sim K'} \smallsetminus K') \enspace,
\end{equation}
where we put $\Pi_{J'}=w_0(J_{\sim K})w_0(J_{\sim K} \smallsetminus K) \cdot \Pi_J$ and $K'=(J \cup K) \smallsetminus J'$.

An element $u \in W$ is called a \emph{right divisor} of $w \in W$ if $\ell(w)=\ell(wu^{-1})+\ell(u)$.
The following lemma is a generalization of \cite[Lemma 4.1]{Bri-How}.
The proof of the original lemma can be easily adapted.
\begin{lem}
[See {\cite[Lemma 4.1]{Bri-How}}]
\label{lem:rightdivisor}
Let $w \in W$ and $J,K \subseteq S$, and suppose that $w \cdot \Pi_J \subseteq \Pi$ and $w \cdot \Pi_K \subseteq \Phi^-$.
Then $J \cap K=\emptyset$, the set $J_{\sim K}$ is of finite type and $w_0(J_{\sim K})w_0(J_{\sim K} \smallsetminus K)$ is a right divisor of $w$.
\end{lem}
A groupoid $G$ played a central role in the paper \cite{Bri-How}.
We recall the definition.
Given $I \subseteq S$, put
\begin{eqnarray*}
\mathscr{J}&=&\{J \subseteq S \mid \Pi_J=w \cdot \Pi_I \mbox{ for some } w \in W\} \enspace;\\
G_{J,K}&=&\{(J,w,K) \mid w \in W \mbox{ and } \Pi_J=w \cdot \Pi_K\}\ \mbox{ for } J,K \in \mathscr{J} \enspace.
\end{eqnarray*}
The multiplication in the groupoid $G=\{G_{J,K}\}_{J,K \in \mathscr{J}}$ is defined by
\begin{displaymath}
(J_1,w,J_2)(J_2,u,J_3)=(J_1,wu,J_3) \enspace.
\end{displaymath}
Note that the map $(I,w,I) \mapsto w$ from $G_{I,I}$ to $N_I$ is a group isomorphism.
Now for $J \subseteq S$ and $s \in S \smallsetminus J$, write
\begin{displaymath}
v[s,J]=w_0(J_{\sim s})w_0(J_{\sim s} \smallsetminus \{s\})\ \mbox{ if } J_{\sim s} \mbox{ is of finite type},
\end{displaymath}
hence $v[s,J] \cdot \Pi_J=\Pi_K$ for a unique $K \subseteq S$.
Moreover, an expression $g=g_1g_2 \cdots g_n$ of $g=(J_0,w,J_n) \in G_{J_0,J_n}$ with $g_i=(J_{i-1},v[s_i,J_i],J_i)$ is called a \emph{standard expression in $G$} if $\ell(w)=\sum_{i=1}^{n}\ell(v[s_i,J_i])$.

The argument in \cite{Bri-How} requires the following theorem of Deodhar:
\begin{thm}
[{\cite[Proposition 5.5]{Deo_root}}]
\label{thm:Deodhar_decomposition}
Any $(J,w,K) \in G$ admits a standard expression in $G$.
Moreover, if $s \in S$ and $w \cdot \alpha_s \in \Phi^-$, then this expression can be chosen in such a way that it ends with $(J',v[s,K],K)$ for some $J'$.
\end{thm}
Now we summarize some results of \cite{Bri-How} required in this paper:
\begin{thm}
\label{thm:resultsofBri-How}
\begin{enumerate}
\item {\rm (\cite[Theorem 2.4]{Bri-How})} Suppose that $s,t \not\in J \subseteq S$, $s \neq t$ and $J'=J_{\sim \{s,t\}}$ is of finite type, hence $w=w_0(J')w_0(J' \smallsetminus \{s,t\})$ maps $\Pi_J$ onto a unique $\Pi_K$.
Then $(K,w,J) \in G$ admits exactly two standard expressions of the form
\begin{displaymath}
(J_0,v[s_1,J_1],J_1) \cdots (J_{n-1},v[s_n,J_n],J_n)
\end{displaymath}
in $G$.
Both of them consist of the same number $n$ of factors and satisfy that $J_i \cup \{s_i\} \subseteq J \cup \{s,t\}$ for all $i$.
Moreover, one of them satisfies $s_n=s$ and the other satisfies $s_n=t$.
\item {\rm (\cite[Theorem A]{Bri-How})} The groupoid $G$ is generated by the elements $(J,v[s,K],K)$.
Moreover, the following two kinds of relations
\begin{itemize}
\item $(J,v[s,K],K)(K,v[t,J],J)=1$, with $K \cup \{s\}=J \cup \{t\}$,
\item $g_1g_2 \cdots g_n=g'_1g'_2 \cdots g'_n$, where $g_1g_2 \cdots g_n$ and $g'_1g'_2 \cdots g'_n$ are the standard expressions of a common $w_0(J')w_0(J' \smallsetminus \{s,t\})$ such that $s,t \not\in J \subseteq S$, $s \neq t$ and $J'=J_{\sim \{s,t\}}$ is of finite type,
\end{itemize}
are fundamental relations of $G$ with respect to these generators.
\end{enumerate}
\end{thm}
Note that $g_1 \neq g'_1$ in the statement 2, since $g_n^{-1} \cdots g_2^{-1}g_1^{-1}$ and $g'_n{}^{-1} \cdots g'_2{}^{-1}g'_1{}^{-1}$ are the standard expressions of the element $(w_0(J')w_0(J' \smallsetminus \{s,t\}))^{-1}$ (see (\ref{eq:inverseofw0w0}) and the statement 1).

In this paper, we say that a generator $(J,v[s,K],K)$ of $G$ is a \emph{loop generator} if $J=K$, namely if it is a loop of the groupoid $G$ (see Section \ref{sec:groupoids} for terminology).
%%%%%
\subsection{Lifting from $G$ to $C$}
\label{sec:groupoidC_lifting}
Now it is straightforward to show that $C$ is a covering groupoid of $G$ with covering map that sends $x \in \mathcal{V}(C)$ to $[x] \in \mathscr{J}$ and $w \in C_{x,y}$ to $([x],w,[y]) \in G_{[x],[y]}$ (see Section \ref{sec:groupoids} for the terminology).
This enables us to lift up the above results on $G$ to its covering groupoid $C$, as in Theorem \ref{thm:presentationofC} below (see Proposition \ref{prop:liftofrelations}).

We prepare notations and terminology.
For a pair $\xi=(x,s)$ of $x \in \mathcal{V}(C)$ and $s \in S \smallsetminus [x]$ such that $[x]_{\sim s}$ is of finite type, we have $y=v[s,[x]] \cdot x \in \mathcal{V}(C)$ (see (\ref{eq:defofwcdoty}) for notation), and the generator $([y],v[s,[x]],[x])$ of $G$ has a unique lift in $C_{y,x}$.
We denote this element in $C_{y,x}$ by $w_x^s$ or $w_\xi$, and write
\begin{displaymath}
\varphi(\xi)=\varphi(x,s)=(y,t) \enspace,
\end{displaymath}
where $t$ is the unique element of $([x] \cup \{s\}) \smallsetminus [y]$.
Theorem \ref{thm:presentationofC}(3) below implies that $\varphi$ is an involutive map from the set of all such pairs $(x,s)$ to itself.
Moreover, we also use the terminology \lq\lq a standard expression of $w$ in $C$'' in a similar way, where the elements $w_x^s$ play the role of the $(J,v[s,K],K)$ in $G$.
Now we have the following results on the groupoid $C$ by applying Proposition \ref{prop:liftofrelations} as mentioned above:
\begin{thm}
\label{thm:presentationofC}
\begin{enumerate}
\item Any $w \in C$ admits a standard expression in $C$.
Moreover, if $s \in S$ and $w \cdot \alpha_s \in \Phi^-$, then this expression can be chosen in such a way that it ends with $w_x^s$.
\item Suppose that $x \in \mathcal{V}(C)$, $s,t \in S \smallsetminus [x]$, $s \neq t$ and $J=[x]_{\sim \{s,t\}}$ is of finite type, hence $w=w_0(J)w_0(J \smallsetminus \{s,t\})$ belongs to $C_{y,x}$ for a unique $y \in \mathcal{V}(C)$.
Then $w$ admits exactly two standard expressions of the form
\begin{displaymath}
w_{z_1}^{s_1}w_{z_2}^{s_2} \cdots w_{z_n}^{s_n}
\end{displaymath}
in $C$.
Both of them consist of the same number $n$ of factors and satisfy that $[z_i] \cup \{s_i\} \subseteq [x] \cup \{s,t\}$ for all $i$.
Moreover, one of them satisfies $s_n=s$ and the other satisfies $s_n=t$.
\item The groupoid $C$ is generated by the elements $w_x^s$.
Moreover, the following two kinds of relations
\begin{itemize}
\item $w_x^sw_y^t=1$, with $[x] \cup \{s\}=[y] \cup \{t\}$ (or equivalently, $w_\xi w_{\varphi(\xi)}=1$),
\item $c_1c_2 \cdots c_n=c'_1c'_2 \cdots c'_n$, where $c_1c_2 \cdots c_n$ and $c'_1c'_2 \cdots c'_n$ are the standard expressions of a common $w_0(J)w_0(J \smallsetminus \{s,t\})$ such that $x \in \mathcal{V}(C)$, $s,t \in S \smallsetminus [x]$, $s \neq t$ and $J=[x]_{\sim \{s,t\}}$ is of finite type,
\end{itemize}
are fundamental relations of $C$ with respect to these generators.
In addition, for the relation of the second type, we have $c_1 \neq c'_1$.
\end{enumerate}
\end{thm}
We use the term \lq\lq loop generator'' also for $C$; an element $w_x^s$ is a loop generator of $C$ if $w_x^s \in C_{x,x}$, or equivalently $w_x^s \cdot x=x$.

We refer to any transformation of expressions in $C$ of the form
\begin{displaymath}
w_1c_1c_2 \cdots c_nw_2 \leadsto w_1c'_1c'_2 \cdots c'_nw_2 \enspace,
\end{displaymath}
where $c_1c_2 \cdots c_n$ and $c'_1c'_2 \cdots c'_n$ are the two expressions in the second relation in Theorem \ref{thm:presentationofC}(3), as a \emph{generalized braid move} (or a \emph{GBM} in short).
Its \emph{loop number} is defined as a half of the total number of loop generators contained in two expressions $c_1c_2 \cdots c_n$ and $c'_1c'_2 \cdots c'_n$.
In fact, this is equal to the number of loop generators in $c_1c_2 \cdots c_n$, or equivalently in $c'_1c'_2 \cdots c'_n$ (see Remark \ref{rem:noincreasing} below).
Now Theorem \ref{thm:presentationofC}(3) says that any two expressions in $C$ of the same element can be converted to each other by generalized braid moves together with insertions and deletions of subwords of the form $w_\xi w_{\varphi(\xi)}$.
This property will be enhanced in Proposition \ref{prop:wordprobleminC}.
%%%%%
\subsection{The graph $\mathcal{C}$}
\label{sec:graphC}
The set of generators $w_x^s$ of $C$ given in Theorem \ref{thm:presentationofC}(3) can be regarded as the edge set of a (connected) graph $\mathcal{C}$ with vertex set $\mathcal{V}(C)$, where $w_x^s$ is an edge from $x \in \mathcal{V}(C)$ to $w_x^s \cdot x \in \mathcal{V}(C)$.
To regard $\mathcal{C}$ as an unoriented graph, we identify each edge $w_\xi$ with its opposite $w_{\varphi(\xi)}$.
Since we write an edge or a path in $\mathcal{C}$ \emph{from right to left} by the convention mentioned in Section \ref{sec:groupoids}, the paths in $\mathcal{C}$ are the expressions of elements in $C$.
This graph $\mathcal{C}$ will be used in our argument below.

\begin{figure}[htbp]
\centering
\begin{picture}(380,45)
\put(15,20){\circle{30}}
\put(30,20){\circle*{10}}
\put(30,20){\line(1,0){30}}
\put(60,20){\circle*{10}}
\put(60,20){\line(1,0){20}}\put(83,17){$\cdots$}\put(100,20){\line(1,0){20}}
\put(120,20){\circle*{10}}
\put(135,20){\circle{30}}
\put(280,15){\line(-1,0){40}}
\put(240,15){\line(0,-1){10}}\put(240,5){\line(-1,0){30}}\put(210,5){\line(0,1){30}}\put(210,35){\line(1,0){30}}\put(240,35){\line(0,-1){10}}
\put(240,25){\line(1,0){100}}
\put(340,25){\line(0,1){10}}\put(340,35){\line(1,0){30}}\put(370,35){\line(0,-1){30}}\put(370,5){\line(-1,0){30}}\put(340,5){\line(0,1){10}}
\put(340,15){\line(-1,0){45}}
\put(280,15){\circle*{10}}
\put(210,15){\hbox to0pt{\hss$\wedge$\hss}}
\put(260,23){\hbox to0pt{\hss$>$\hss}}
\put(370,15){\hbox to0pt{\hss$\vee$\hss}}
\put(297,13){\hbox to0pt{\hss$<$\hss}}
\end{picture}
\caption{The component $\mathcal{G}$ of the second type}
\label{fig:shapeofcomponent}
\end{figure}
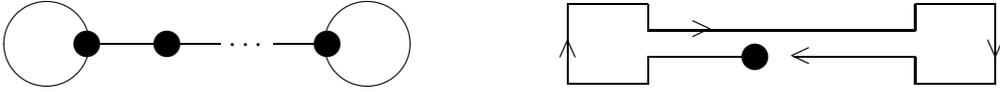

Let $x \in \mathcal{V}(C)$ and $s,t \in S \smallsetminus [x]$ such that $s \neq t$ and $J_{\sim (J \smallsetminus [x])}$ is of finite type, where $J=[x] \cup \{s,t\}$.
Let $\mathcal{C}(J)$ be the subgraph of $\mathcal{C}$ consisting of all vertices $y \in \mathcal{V}(C)$ with $[y] \subseteq J$ and all edges $w_y^{s'}$ with $[y] \cup \{s'\} \subseteq J$.
Moreover, let $\mathcal{G}$ be the connected component of $\mathcal{C}(J)$ containing $x$.
Then it is easy to show that, for any vertex $y$ of $\mathcal{G}$, we have $|J \smallsetminus [y]|=2$, $J_{\sim (J \smallsetminus [y])}=J_{\sim (J \smallsetminus [x])}$, and $y_\lambda=x_\lambda$ whenever $\lambda \in \Lambda$ and $x_\lambda \not\in J_{\sim (J \smallsetminus [x])}$.
(Indeed, it suffices to check the property only for adjacent vertices $x,y$ in $\mathcal{G}$.)
Thus $\mathcal{G}$ is a finite graph in which every vertex is adjacent to exactly two edges, hence it is classified into the following two types:
\begin{enumerate}
\item A cycle without loops.
We refer to any nonempty, non-backtracking closed path $p$ in this $\mathcal{G}$ as a \emph{circular tour}.
\item A union of two loops and a simple (possibly empty) path joining the loops (see the left-hand side of Figure \ref{fig:shapeofcomponent}).
We say that a closed path $p$ in this $\mathcal{G}$ is a \emph{shuttling tour} if it visits every vertex of $\mathcal{G}$ exactly twice (except the start point of the path) and passes each of the two loops exactly once (see the right-hand side of Figure \ref{fig:shapeofcomponent}).
\end{enumerate}
In both cases, the path $p$ represents an element of the finite parabolic subgroup $W_{J_{\sim (J \smallsetminus [x])}}$, hence the order of $p$ as an element of $W$ is finite.
%%%%%
%%%%%

The next lemma relates the shuttling tours with the generalized braid moves:
\begin{lem}
\label{lem:shuttlingandGBM}
Let $q$ be a path in $\mathcal{C}$ and $k \geq 1$.
Then the followings are equivalent:
\begin{enumerate}
\item $q=p^k$ as paths in $\mathcal{C}$ for a shuttling tour $p$ of order $k$;
\item $q=(c'_1 \cdots c'_r)^{-1}c_1 \cdots c_r$ as paths in $\mathcal{C}$ for a generalized braid move $c_1 \cdots c_r \leadsto c'_1 \cdots c'_r$ of loop number $k$.
\end{enumerate}
\end{lem}
\begin{proof}
Assuming the property 1, let $p$ correspond to the graph $\mathcal{G}$, start with $w_y^{s'}$ and end with $(w_y^{t'})^{-1}$.
Then $[y]_{\sim \{s',t'\}}$ is of finite type as mentioned above, therefore the triple $(y,s',t')$ yields a GBM $c_1 \cdots c_r \leadsto c'_1 \cdots c'_r$ such that $(c_r,c'_r)=(w_y^{s'},w_y^{t'})$, both $c_1 \cdots c_r$ and $c'_1 \cdots c'_r$ are non-backtracking distinct paths in $\mathcal{G}$, and $c_1 \neq c'_1$.
Now the shape of $\mathcal{G}$ forces the closed path $q'=(c'_1 \cdots c'_r)^{-1}c_1 \cdots c_r$ to be a power of $p$, say, $p^n$ with $n>0$.
Then this GBM has loop number $n$, since the path $p^n$ contains $2n$ loops.

We show that $n=k$.
Since $q'=p^n$ represents an identity element in $C$, the order $k$ of $p$ is a divisor of $n$.
Now if $k<n$, then $k \leq n/2$ and a path $p^k$ representing an identity element must be a subpath of the first half $c_1 \cdots c_r$ of $q'=p^n$, contradicting the standardness of $c_1 \cdots c_r$.
Thus $n=k$ as desired, hence the property 2 follows.

On the other hand, assuming the property 2, let the triple $(x,s,t)$ correspond to the GBM as in Theorem \ref{thm:presentationofC}(3).
Then, since the closed path $q=(c'_1 \cdots c'_r)^{-1}c_1 \cdots c_r$ contains $2k$ loops, the same argument as the first paragraph implies that $q$ is a certain power $p^n$ of a shuttling tour $p$ with $n>0$, this $p$ has order $n$, and we have $n=k$ since $p^n$ contains $2n$ loops.
Hence the property 1 follows, concluding the proof of Lemma \ref{lem:shuttlingandGBM}.
\end{proof}
%%
%%%%%
\section{The decomposition of $C$}
\label{sec:decompofC}
The aim of this section is to prove the decomposition (\ref{eq:decompofC_I}) of the group $C_I$ and describe the factors in detail.
In the course of our argument (see Section \ref{sec:analogyofC}), we develop more similarities of $C$ with Coxeter groups, in addition to those inherited from those of $G$ shown in \cite{Bri-How} and summarized above.
%%%%%
\subsection{Similarities of $C$ with Coxeter groups}
\label{sec:analogyofC}
We start with the following lemma, which says that any loop generator of $C$ lying in $C_{x,x}$ is a reflection along a root in $\Phi^{\perp [x]}$.
This property of our groupoid $C$, which the groupoid $G$ does not possess, plays an important role in our argument below.
\begin{lem}
\label{lem:charofBphi}
For $x \in \mathcal{V}(C)$ and $s \in S \smallsetminus [x]$, the three conditions are equivalent:
\begin{enumerate}
\item $[x]_{\sim s}$ is of finite type, and $\varphi(x,s)=(x,s)$, i.e., $w_x^s \in C_{x,x}$;
\item $[x]_{\sim s}$ is of finite type, and $\Phi^{\perp [x]}[w_x^s] \neq \emptyset$, i.e., $w_x^s \not\in Y$;
\item $\Phi_{[x] \cup \{s\}}^{\perp [x]} \neq \emptyset$.
\end{enumerate}
If these conditions are satisfied, we have $\Phi^{\perp [x]}[w_x^s]=(\Phi_{[x] \cup \{s\}}^{\perp [x]})^+=\{\gamma(x,s)\}$ for a unique positive root $\gamma(x,s)$ such that $s_{\gamma(x,s)}=w_x^s$.
\end{lem}
\begin{proof}
Put $\xi=(x,s)$.
First, we deduce the property 2 from the property 1.
The property 1 implies that $w_\xi=w_{\varphi(\xi)}$, while $w_{\varphi(\xi)}=w_\xi^{-1}$ (see Theorem \ref{thm:presentationofC}(3)), therefore $w_\xi{}^2=1$.
Thus Theorem \ref{thm:involution} gives us a root $\gamma$ with $w_\xi \cdot \gamma=-\gamma$, which satisfies that $\gamma \in \Phi^{\perp [x]}$ since $w_\xi \in C_{x,x}$ and the form $\langle \,,\, \rangle$ is $W$-invariant, proving the property 2 as desired.
Moreover, the property 2 implies the property 3 since $\Phi[w_\xi] \subseteq \Phi_{[x] \cup \{s\}}$.
This inclusion also implies that $\Phi^{\perp [x]}[w_\xi] \subseteq (\Phi_{[x] \cup \{s\}}^{\perp [x]})^+$.

Now we deduce the property 1 from the property 3.
Let $\gamma \in (\Phi_{[x] \cup \{s\}}^{\perp [x]})^+$.
Then we have $s_\gamma \in C_{x,x} \cap W_{[x] \cup \{s\}}$, therefore $\emptyset \neq \Phi[s_\gamma] \subseteq \Phi_{[x] \cup \{s\}}^+ \smallsetminus \Phi_{[x]}$, forcing the root $s_\gamma \cdot \alpha_s$ to be negative.
Hence Theorem \ref{thm:presentationofC}(1) gives us a standard expression of $s_\gamma$ in $C$ ending with $w_\xi$.
Moreover, this expression consists of only the term $w_\xi$, since the equality $\Phi_{[x] \cup \{s\}}^+ \smallsetminus \Phi_{[x]}=\Phi[w_\xi]$ implies that $\Phi[s_\gamma] \subseteq \Phi[w_\xi]$ and $\ell(s_\gamma) \leq \ell(w_\xi)$.
This proves that $s_\gamma=w_\xi \in C_{x,x}$.
Since $\gamma \in (\Phi_{[x] \cup \{s\}}^{\perp [x]})^+$ was chosen arbitrarily, the uniqueness of $\gamma=\gamma(x,s)$ and the inclusion $(\Phi_{[x] \cup \{s\}}^{\perp [x]})^+ = \{\gamma(x,s)\} \subseteq \Phi^{\perp [x]}[w_\xi]$ follow.
Moreover, it also follows that $w_\xi \cdot x=x$, therefore $\varphi(\xi)=\xi$ as desired.
Hence the proof of Lemma \ref{lem:charofBphi} is concluded.
\end{proof}
The next property is an analogy of the Exchange Condition for Coxeter groups:
\begin{lem}
\label{lem:weakExchangeConditionforC}
Let $w=w_{\xi_1}w_{\xi_2} \cdots w_{\xi_n}w_\zeta$ be an expression in $C$, where $\zeta=(x,s)$ and $\varphi(\zeta)=(y,t)$, such that the expression $w_{\xi_1} \cdots w_{\xi_n}$ is standard.
Then the conditions 1--3 below are equivalent in general, and moreover, the condition 4 is also equivalent to the first three when $w_\zeta=s_\gamma$ is a loop generator with $\gamma$ a positive root (see Lemma \ref{lem:charofBphi}):
\begin{enumerate}
\item the expression $w_{\xi_1} \cdots w_{\xi_n}w_\zeta$ is not standard;
\item $w_{\xi_1} \cdots w_{\xi_n} \cdot \alpha_t \in \Phi^-$;
\item $w_{\xi_1} \cdots w_{\xi_n} \in C$ admits a standard expression ending with $w_\zeta^{-1}=w_{\varphi(\zeta)}$;
\item $w_{\xi_1} \cdots w_{\xi_n} \cdot \gamma \in \Phi^-$.
\end{enumerate}
\end{lem}
\begin{proof}
Put $u=w_{\xi_1} \cdots w_{\xi_n}$.
First, the condition 2 implies the condition 3 by Theorem \ref{thm:presentationofC}(1), while the condition 3 implies the condition 1 since now $\ell(uw_\zeta)=\ell(u)-\ell(w_\zeta)$.
We show that the condition 1 implies the condition 2.
By the hypothesis and the condition 1, we have $\ell(uw_\zeta)<\ell(u)+\ell(w_\zeta)$, hence $u$ maps some $\gamma \in \Phi[w_\zeta^{-1}]=\Phi_{[y] \cup \{t\}}^+ \smallsetminus \Phi_{[y]}$ to a negative root (see Lemma \ref{lem:lengthofmultiple}).
Since $u$ maps $\Pi_{[y]}$ into $\Pi$, it must map $\alpha_t$ to a negative root, proving the condition 2 as desired.

In the special case $w_\zeta=s_\gamma$, the condition 4 implies the condition 1, since now both $u$ and $s_\gamma$ send $\gamma \in \Phi^+$ into $\Phi^-$ and hence $\ell(us_\gamma)<\ell(u)+\ell(s_\gamma)$ by Lemma \ref{lem:lengthofmultiple}.
On the other hand, the condition 3 implies the condition 4 since now we have $\gamma \in \Phi[w_\zeta^{-1}] \subseteq \Phi[u]$ by Lemma \ref{lem:lengthofmultiple} again.
Hence the proof of Lemma \ref{lem:weakExchangeConditionforC} is concluded.
\end{proof}
Given $w \in C$ and its standard expression, let $\mathsf{lp}(w)$ denote the number of the loop generators contained in this expression.
It will be shown in Proposition \ref{prop:looplengthofw} that $\mathsf{lp}(w)$ is well-defined regardless of the choice of the standard expression of $w$.
The next property is an analogy of a well-known property of Coxeter groups concerning the length of $w$ and the set $\Phi[w]$.
The proof is also analogous to the case of Coxeter groups.
\begin{prop}
\label{prop:looplengthofw}
Let $w \in C_{y,x}$.
Then $\mathsf{lp}(w)$ is equal to the cardinality of the set $\Phi^{\perp [x]}[w]$, hence is determined just by $w$ regardless of the given standard expression.
If $w$ admits a standard expression of the form $u_0s_{\gamma_1}u_1s_{\gamma_2}u_2 \cdots u_{n-1}s_{\gamma_n}u_n$, where each $u_i$ contains no loop generators and each $s_{\gamma_i}$ is a loop generator with $\gamma_i$ positive, then
\begin{displaymath}
\Phi^{\perp [x]}[w]=\{\beta_i=(u_is_{\gamma_{i+1}}u_{i+1} \cdots s_{\gamma_n}u_n)^{-1} \cdot \gamma_i \mid 1 \leq i \leq n\} \enspace.
\end{displaymath}
\end{prop}
\begin{proof}
Put $w_{i,j}=u_is_{\gamma_{i+1}}u_{i+1} \cdots s_{\gamma_j}u_j$ for indices $i \leq j$, and $z_i=w_{i,n} \cdot x$.
First we show that all the $\beta_i$ in the statement are distinct and lie in $\Phi^{\perp [x]}[w]$, proving that $n \leq |\Phi^{\perp [x]}[w]|$.
Note that $\beta_i \in \Phi^{\perp [x]}$ since $\gamma_i \in \Phi^{\perp [z_i]}$.
Now if $1 \leq i<j \leq n$ and $\beta_i=\beta_j$, then we have
\begin{displaymath}
w_{i,j-1}{}^{-1} \cdot \gamma_i=s_{\gamma_j}w_{j,n} \cdot \beta_i=s_{\gamma_j}w_{j,n} \cdot \beta_j = s_{\gamma_j} \cdot \gamma_j = -\gamma_j \in \Phi^- \enspace,
\end{displaymath}
therefore Lemma \ref{lem:weakExchangeConditionforC} gives us a standard expression of $w_{i,j-1}$, hence of $w_{i,n}$, beginning with $s_{\gamma_i}$.
The same situation occurs if $\beta_i \in \Phi^-$.
However, this prevents the given expression of $w$ from being standard, contradicting the hypothesis.
Thus all $\beta_i$ are positive and distinct.
Moreover, Lemma \ref{lem:lengthofmultiple} implies that $\beta_i \in \Phi[s_{\gamma_i}w_{i,n}] \subseteq \Phi[w]$.
This completes the first claim of this proof.

From now, we prove the other inequality $|\Phi^{\perp [x]}[w]| \leq n$ by induction on $n$.
First, Lemma \ref{lem:charofBphi} says that any non-loop generator of $C$, hence the $u_n$, lies in the groupoid $Y$.
This proves the claim for $n=0$ and allows us to assume that $u_n=1$ without loss of generality.
Now since $\Phi^{\perp [x]}[s_{\gamma_n}]=\{\gamma_n\}$ (see Lemma \ref{lem:charofBphi}), $s_{\gamma_n}$ maps $\Phi^{\perp [x]}[w] \smallsetminus \{\gamma_n\}$ into the set $\Phi^{\perp [x]}[w_{0,n-1}]$ with cardinality not larger than $n-1$ (the induction assumption), proving the desired inequality.
Hence the proof of Proposition \ref{prop:looplengthofw} is concluded.
\end{proof}
The next result shows a further remarkable property of $C$, enhancing Theorem \ref{thm:presentationofC}(3).
Note that the corresponding property of the fundamental relations of Coxeter groups plays a crucial role in the solution of the word problem in Coxeter groups; see e.g., \cite[Theorem 3.3.1]{Bjo-Bre}.
The proof is analogous to the case of Coxeter groups again.
\begin{prop}
\label{prop:wordprobleminC}
Let $w \in C$ be any element.
\begin{enumerate}
\item Any two standard expressions of $w$ can be converted to each other by using the generalized braid moves only.
\item Any expression of $w$ can be converted to a given standard expression of $w$ by using the generalized braid moves and cancellations of subwords of the form $w_\xi w_{\varphi(\xi)}$ only, not using insertions of subwords $w_{\xi}w_{\varphi(\xi)}$.
\end{enumerate}
\end{prop}
\begin{proof}
For the claim 1, we proceed the proof by induction on $\ell(w)$, the case $\ell(w)=0$ being trivial.
Suppose that $\ell(w)>0$, and let $w_{x_1}^{s_1} \cdots w_{x_n}^{s_n}$ and $w_{y_1}^{t_1} \cdots w_{y_m}^{t_m}$ be two standard expressions of $w$, hence $x_n=y_m$.
Then we have $w \cdot \alpha_{s_n} \in \Phi^-$ by the fact $w_{x_n}^{s_n} \cdot \alpha_{s_n} \in \Phi^-$ and Lemma \ref{lem:lengthofmultiple}, and $w \cdot \alpha_{t_m} \in \Phi^-$ similarly.
Now the combination of Lemma \ref{lem:rightdivisor} and Theorem \ref{thm:presentationofC}(1)(2) implies that $J=[x_n]_{\sim \{s_n,t_m\}}$ is of finite type and there are two standard expressions of $w$ of the forms $uc_1 \cdots c_r$ and $uc'_1 \cdots c'_r$, where $c_1 \cdots c_r$ and $c'_1 \cdots c'_r$ are the standard expressions of $w_0(J)w_0(J \smallsetminus \{s_n,t_m\}) \in C$ ending with $c_r=w_{x_n}^{s_n}$ and $c'_r=w_{y_m}^{t_m}$, respectively.
Then we have
\begin{displaymath}
w_{x_1}^{s_1} \cdots w_{x_n}^{s_n} \overset{GBMs}{\leadsto} uc_1 \cdots c_r \overset{GBM}{\leadsto} uc'_1 \cdots c'_r \overset{GBMs}{\leadsto} w_{y_1}^{t_1} \cdots w_{y_m}^{t_m}
\end{displaymath}
as desired, where the first and the third transformations come from the induction assumption and the second one is the definition of the GBM.

For the claim 2, let $w_{\xi_1} \cdots w_{\xi_n}$ be an expression of $w$.
We show the claim by induction on the total length of the $w_{\xi_i}$, the case of length zero being trivial.
Owing to the claim 1, we may assume that this expression is not standard.
Take the last index $i$ such that $w_{\xi_1} \cdots w_{\xi_{i-1}}$ is standard.
Then by Lemma \ref{lem:weakExchangeConditionforC}, $w_{\xi_1} \cdots w_{\xi_{i-1}}$ admits another standard expression of the form $w_{\zeta_1} \cdots w_{\zeta_m}w_{\varphi(\xi_i)}$, which can be reached from $w_{\xi_1} \cdots w_{\xi_{i-1}}$ by GBMs only (apply the claim 1).
Thus we have
\begin{displaymath}
w_{\xi_1} \cdots w_{\xi_n} \overset{GBMs}{\leadsto} w_{\zeta_1} \cdots w_{\zeta_m}w_{\varphi(\xi_i)}w_{\xi_i} \cdots w_{\xi_n} \overset{cancel}{\leadsto} w_{\zeta_1} \cdots w_{\zeta_m}w_{\xi_{i+1}} \cdots w_{\xi_n}
\end{displaymath}
and the total length decreases through the transformation.
Hence the induction works and the claim follows.
\end{proof}
\begin{rem}
\label{rem:noincreasing}
As shown in Theorem \ref{thm:presentationofC}(2), both terms of a generalized braid move have the same number of generators, the same number of loop generators (see Proposition \ref{prop:looplengthofw}) and the same total length of generators contained.
Thus none of those quantities increases in a transformation appearing in Proposition \ref{prop:wordprobleminC}.
As a result, any expression of $w \in C$ contains at least $\mathsf{lp}(w)$ loop generators.
\end{rem}
%%
%%%%%
\subsection{The factorization of $C$}
\label{sec:decompofC_factors}
The previous results enable us to deduce the following properties.
Recall from Section \ref{sec:groupoidC_def} that $W^{\perp [x]}$ is normal in $C_{x,x}$ and $W^{\perp [x]} \cap Y_{x,x}=1$.
\begin{thm}
\label{thm:decompofC}
\begin{enumerate}
\item We have
\begin{displaymath}
C_{y,x}=W^{\perp [y]} \cdot Y_{y,x}=Y_{y,x} \cdot W^{\perp [x]} \mbox{ for all } x,y \in \mathcal{V}(C) \enspace.
\end{displaymath}
Hence $C_{x,x}$ admits a semidirect product decomposition $C_{x,x}=W^{\perp [x]} \rtimes Y_{x,x}$.
\item The groupoid $Y$ is generated by all the non-loop generators of $C$.
We have
\begin{displaymath}
Y_{y,x}=\{w \in C_{y,x} \mid \mathsf{lp}(w)=0\} \mbox{ for all } x,y \in \mathcal{V}(C) \enspace.
\end{displaymath}
Moreover, Proposition \ref{prop:wordprobleminC} holds also for $Y$ under the modification that the generalized braid moves are restricted to the ones containing no loop generators.
\item The length of an element $w \in W^{\perp [x]}$ with respect to the generating set $R^x = S(\Phi^{\perp [x]})$ of $W^{\perp x}$ is equal to $\mathsf{lp}(w)$, and we have
\begin{equation}
\label{eq:elementsofRx}
R^{[x]}=\{ws_\gamma w^{-1} \mid s_\gamma \in C_{y,y} \mbox{ and } w \in Y_{x,y} \mbox{ for some } y \in \mathcal{V}(C)\} \enspace.
\end{equation}
Moreover, $w \in W^{\perp [x]}$ admits an expression $(w_1s_{\gamma_1}{w_1}^{-1}) \cdots (w_ns_{\gamma_n}{w_n}^{-1})$ with $n=\mathsf{lp}(w)$ such that the expression obtained from it by replacing every $w_1$, ${w_i}^{-1}w_{i+1}$, ${w_n}^{-1}$ with its standard expression is also standard.
\item If $w \in Y_{y,x}$, then the map $u \mapsto wuw^{-1}$ is an isomorphism of Coxeter systems from $(W^{\perp [x]},R^{[x]})$ to $(W^{\perp [y]},R^{[y]})$.
\end{enumerate}
\end{thm}
\begin{proof}
First we note the following fact used in the proof.
If $w \in C$ admits an expression of the form $u_0s_{\gamma_1}u_1s_{\gamma_2} \cdots u_{n-1}s_{\gamma_n}u_n$, where the $s_{\gamma_i}$ are loop generators and the $u_i$ contain no loop generators, then
\begin{equation}
\label{eq:decompofC_1}
w=(w_1s_{\gamma_1}{w_1}^{-1} \cdots w_ns_{\gamma_n}{w_n}^{-1})w_{n+1}=(s_{w_1 \cdot \gamma_1} \cdots s_{w_n \cdot \gamma_n})w_{n+1}
\end{equation}
where we put $w_i=u_0u_1 \cdots u_{i-1}$ for each $1 \leq i \leq n+1$.

For the claim 1, it suffices to show the first part of the claim, in particular the first equality, since then the second equality follows by taking the inverse.
This is done by using the transformation (\ref{eq:decompofC_1}) and Lemma \ref{lem:charofBphi} (indeed, if $w \in C_{y,x}$, then $w_i \cdot \gamma_i \in \Phi^{\perp [y]}$ and $w_{n+1} \in Y_{y,x}$).

For the claim 2, the first and the third parts follow from the second part, Theorem \ref{thm:presentationofC}(1) and Remark \ref{rem:noincreasing}.
Moreover, the second one is deduced from Proposition \ref{prop:looplengthofw}.

For the claim 3, note that the claim on the length of $w$ follows from Theorem \ref{thm:reflectionsubgroup_Deodhar}(5) and Proposition \ref{prop:looplengthofw}.
Moreover, any element $u$ of the right-hand side of (\ref{eq:elementsofRx}) satisfies that $\mathsf{lp}(u) \leq 1$ (see Remark \ref{rem:noincreasing}) and $1 \neq u \in W^{\perp [x]}$, hence $\mathsf{lp}(u)=1$ and $u \in R^{[x]}$ by the above claim.
Now take a standard expression $u_0s_{\gamma_1}u_1s_{\gamma_2} \cdots u_{n-1}s_{\gamma_n}u_n$ of $w \in W^{\perp [x]}$ in $C$ of the above form, hence $n=\mathsf{lp}(w)$.
Then in the equality (\ref{eq:decompofC_1}), we have $w_{n+1}=1$ since $W^{\perp [x]} \cap Y_{x,x}=1$, while $w_1=u_0$ and ${w_i}^{-1}w_{i+1}=u_i$ for all $1 \leq i \leq n$.
Hence (\ref{eq:decompofC_1}) is the desired expression of $w$.
This argument also implies that the set $R^{[x]}=\{w \in W^{\perp [x]} \mid \mathsf{lp}(w)=1\}$ is contained in the right-hand side of (\ref{eq:elementsofRx}).
Hence the claim holds.

Finally, the claim 4 is a corollary of the claim 3.
Hence Theorem \ref{thm:decompofC} holds.
\end{proof}
\begin{rem}
\label{rem:relationofCandG}
It is shown in \cite{Bri-How} that the vertex group $G_I = G_{I,I}$ of the groupoid $G$ admits a similar decomposition $G_I=\widetilde{N}_I \rtimes M_I$ into a Coxeter group $\widetilde{N}_I$ and a vertex group $M_I$ of a groupoid $M$, where $M$ and $\widetilde{N}_I$ are (by definition) generated by non-loop generators of $G$, and by conjugates of loop generators of $G$ by elements of $M$, respectively.
However, in contrast with the case of $C$, the properties of $W^{\perp I}$ and $Y$ are not inherited immediately from $\widetilde{N}_I$ and $M$, since a lift of a loop generator of $G$ may be a non-loop generator of $C$ and not all the generators of $W^{\perp I}$ are lifts of those of $\widetilde{N}_I$.
More precisely, it can be shown that $W^{\perp I}$ is a normal reflection subgroup of the Coxeter group $\widetilde{N}_I$, but \emph{not} a parabolic one; hence the structure of $W^{\perp I}$ is still hard to describe even if that of $\widetilde{N}_I$ is clear.
\end{rem}
%%
%%%%%
\subsection{The factor $Y$}
\label{sec:decompofC_factorY}
In this and the following subsections, we study the structure of the two factors $W^{\perp I}$ and $Y_I$ more precisely.
First, we have the following result:
\begin{prop}
\label{prop:Yistorsionfree}
Each vertex group $Y_{x,x}$ of $Y$ is torsion-free.
\end{prop}
\begin{proof}
If $w \in Y_{x,x}$ has finite order, then Proposition \ref{prop:fixingroots} gives us a decomposition $w=s_{\beta_1} \cdots s_{\beta_r}$ such that each $s_{\beta_j}$ fixes $\Pi_{[x]}$ pointwise, therefore we have $\beta_j \in \Phi^{\perp [x]}$ and $w \in Y_{x,x} \cap W^{\perp [x]}=1$, implying $w = 1$ as desired.
\end{proof}
This yields the following analogy of Lemma \ref{lem:shuttlingandGBM} (see Section \ref{sec:graphC} for the definitions):
\begin{cor}
\label{cor:circulartourandGBM}
Any circular tour has order $1$.
Moreover, for a path $p$ in the graph $\mathcal{C}$, the followings are equivalent:
\begin{enumerate}
\item $p$ is a circular tour;
\item $p=(c'_1 \cdots c'_r)^{-1}c_1 \cdots c_r$ as paths in $\mathcal{C}$ for a generalized braid move $c_1 \cdots c_r \leadsto c'_1 \cdots c'_r$ with loop number $0$.
\end{enumerate}
\end{cor}
\begin{proof}
By definition, any circular tour $p$ represents an element of finite order in some $Y_{x,x}$ that is a torsion-free group (see Proposition \ref{prop:Yistorsionfree}), hence the element must be identity.
The second claim is deduced by the same argument as Lemma \ref{lem:shuttlingandGBM}.
\end{proof}
The combination of this corollary and Theorem \ref{thm:decompofC}(2) yields the next description of the groupoid $Y$.
Define a complex $\mathcal{Y}$ (see Section \ref{sec:groupoids} for terminology) such that its $1$-skeleton $\mathcal{Y}^1$ is the subgraph of $\mathcal{C}$ obtained by deleting all the loops, and $\mathcal{Y}$ has a $2$-cell with boundary $\mathcal{G}$ for each cycle $\mathcal{G}$ as in the definition of circular tours (see Section \ref{sec:graphC}).
Note that every path in $\mathcal{Y}$ represents an element of $C$.
\begin{thm}
\label{thm:presentationofY}
The groupoid $Y$ is naturally isomorphic to the fundamental groupoid $\pi_1(\mathcal{Y};*,*)$ of the complex $\mathcal{Y}$.
Hence each $Y_{x,x}$ is also isomorphic to the fundamental group $\pi_1(\mathcal{Y};x)$ of $\mathcal{Y}$ at $x$.
\end{thm}
Note that the groupoid $C$ also admits a similar description in terms of a certain complex with $\mathcal{C}$ being the $1$-skeleton, though $\mathcal{C}$ involves a loop in general and the boundary of a $2$-cell is not always a simple closed path any longer.

For a further description, we fix a maximal tree $\mathcal{T}$ in the connected graph $\mathcal{Y}^1$.
For $y,z \in \mathcal{V}(C)$, let $p_{z,y}$ denote the unique non-backtracking path in $\mathcal{T}$ from $y$ to $z$.
Note that $p_{z,y}p_{y,x}=p_{z,x}$ and $p_{y,z}=p_{z,y}{}^{-1}$ hold in $\pi_1(\mathcal{Y}^1;*,*)$.
Moreover, put
\begin{equation}
\label{eq:extensionofpath}
q_{(x)}=p_{x,z}qp_{y,x}\ \mbox{ for } x \in \mathcal{V}(C) \mbox{ and a path } q \mbox{ in } \mathcal{Y}^1 \mbox{ from } y \mbox{ to } z \enspace,
\end{equation}
which is the extension of the path $q$ to a vertex $x$ along the tree $\mathcal{T}$.
Note that $(e_{\xi_n} \cdots e_{\xi_1})_{(x)}=(e_{\xi_n})_{(x)} \cdots (e_{\xi_1})_{(x)}$ in $\pi_1(\mathcal{Y}^1;x)$.
Then a theorem in combinatorial group theory (see e.g., \cite[Theorem 5.17]{Coh}) yields the following presentation of $\pi_1(\mathcal{Y};x)$:
\begin{thm}
\label{thm:presentationofpi1Y}
The group $\pi_1(\mathcal{Y};x)$ admits a presentation with generators given by
\begin{displaymath}
(w_y^s)_{(x)}=p_{x,z}w_y^sp_{y,x}\,, \mbox{ with } w_y^s \in E(\mathcal{Y}) \mbox{ and } z=w_y^s \cdot y \enspace,
\end{displaymath}
where $E(\mathcal{Y})$ denotes the set of the oriented edges of $\mathcal{Y}$, and fundamental relations given by
\begin{itemize}
\item $(w_{\xi})_{(x)}(w_{\varphi(\xi)})_{(x)}=1$ for every $w_\xi \in E(\mathcal{Y})$;
\item $(w_{\xi})_{(x)}=1$ for every $w_{\xi} \in E(\mathcal{Y})$ contained in the tree $\mathcal{T}$;
\item $(e_{\xi_n})_{(x)} \cdots (e_{\xi_1})_{(x)}=1$ for a boundary $e_{\xi_n} \cdots e_{\xi_1}$ of each $2$-cell of $\mathcal{Y}$.
\end{itemize}
\end{thm}
\begin{exmp}
\label{ex:caseofsinglereflection}
Here we consider the case $|I|=1$ examined by Brink \cite{Bri}; hence $S^{(\Lambda)} = S$ and $I=\{x_I\}$.
Then for $x \in S$ and $s \in S \smallsetminus \{x\}$, the set $\{x\}_{\sim s}$ is of finite type and $w_x^s \cdot x \neq x$ if and only if $m_{x,s}$ is odd.
Thus the $1$-skeleton $\mathcal{Y}^1$ is the connected component of the \emph{odd Coxeter graph} $\Gamma^{\mathrm{odd}}$ of $(W,S)$ containing $x_I$ used in Brink's description of $Z_W(x_I)$, where $\Gamma^{\mathrm{odd}}$ is (by definition) the subgraph of $\Gamma$ obtained by removing the edges with non-odd labels.
Moreover, the acyclicness of Coxeter graphs of finite type implies that this $\mathcal{Y}$ has no $2$-cell and $\mathcal{Y}=\mathcal{Y}^1$.
Hence our result agrees with the result of \cite{Bri} for this special case.
\end{exmp}
%%
%%%%%
\subsection{The factor $W^{\perp I}$}
\label{sec:decompofC_factorWperpI}
First, for $x \in \mathcal{V}(C)$, define
\begin{displaymath}
\mathcal{R}_x=\{(w,u) \mid w \in Y_{x,x} \mbox{ and } u \mbox{ is a loop generator of } C\}
\end{displaymath}
and
\begin{displaymath}
r_x(w,u)=wu_{(x)}w^{-1} \mbox{ for } (w,u) \in \mathcal{R}_x \mbox{ (see (\ref{eq:extensionofpath}) for the notation)},
\end{displaymath}
hence the generating set $R^{[x]}$ of $W^{\perp [x]}$ consists of all the $r_x(w,u)$ with $(w,u) \in \mathcal{R}_x$ (see Theorem \ref{thm:decompofC}(3)).
The subscripts \lq $x$' are omitted when $x=x_I$.
The pair $(w,u)$ is also denoted by $(w,\gamma)$ when $u=s_\gamma$.
To study relations between the generators $r_x(w,u)$, for integer $1 \leq k<\infty$, we define a relation $\overset{k}{\sim}$ on $\mathcal{R}_x$ by
\begin{displaymath}
(w,\gamma) \overset{k}{\sim} (wq_{(x)},\beta) \mbox{ if there is a shuttling tour of order $k$ of the form } s_\gamma qs_\beta q^{-1}
\end{displaymath}
(see Section \ref{sec:graphC} for the terminology).
Note that the relation $\overset{k}{\sim}$ is symmetric, and the product $r_x(w,\gamma)r_x(u,\beta)$ of two generators of $W^{\perp [x]}$ has order $k$ if $(w,\gamma) \overset{k}{\sim} (u,\beta)$.
Let $\sim$ denote the transitive closure of $\overset{1}{\sim}$, which is an equivalence relation.

The next result describes the structure of the group $W^{\perp [x]}$ in terms of those relations.
Note that the first claim of this result can be deduced directly from Theorem \ref{thm:presentationofC}(3) in a similar way, without Proposition \ref{prop:wordprobleminC} that is essential in the proof of the second claim.
\begin{thm}
\label{thm:presentationofWperp}
The group $W^{\perp [x]}$ admits a presentation with generators given by
\begin{displaymath}
r_x(w,\gamma)\,, \mbox{ with } (w,\gamma) \in \mathcal{R}_x
\end{displaymath}
and fundamental relations given by
\begin{itemize}
\item $r_x(w,\gamma)^2=1$ for every $(w,\gamma) \in \mathcal{R}_x$;
\item $(r_x(w,\gamma)r_x(u,\beta))^k=1$ for every pair $(w,\gamma) \overset{k}{\sim} (u,\beta)$ with $1 \leq k<\infty$.
\end{itemize}
Moreover, any expression of an element of $W^{\perp [x]}$ (with respect to these generators) can be converted into an arbitrarily given shortest expression of the element (with respect to these generators) by using the following two kinds of transformations:
\begin{itemize}
\item $r_x(w,\gamma)^2 \leadsto 1$, corresponding to the first relation above;
\item $r_x(w,\gamma)r_x(u,\beta)r_x(w,\gamma) \cdots \,\leadsto\, r_x(u,\beta)r_x(w,\gamma)r_x(u,\beta) \cdots$ (where both terms consist of $k$ generators), corresponding to the second relation above.
\end{itemize}
\end{thm}
\begin{proof}
Owing to Theorem \ref{thm:decompofC}(3), it suffices to prove that any expression of an element of $W^{\perp [x]}$ can be converted into a given shortest expression \emph{of the form given in Theorem \ref{thm:decompofC}(3)} in the above way, since the above transformations of second type are invertible.
In this proof, we write $w \equiv u$ if two expressions $w$ and $u$ represent the same element, distinguishing from the case $w=u$ where the expressions $w$ and $u$ themselves are equal.

To prove the above claim, write the given original expression as
\begin{equation}
\label{eq:presentationofWperp_1}
w_1l_1w_1^{-1}w_2l_2w_2^{-1} \cdots w_nl_nw_n^{-1}, \mbox{ where } (w_i,l_i) \in \mathcal{R}_x \enspace.
\end{equation}
Put $w_0=1$ and $w_{n+1}=1$.
Then, owing to the assumption on the property of the target shortest expression, the combination of Proposition \ref{prop:wordprobleminC} and Theorem \ref{thm:decompofC}(2) implies the following fact: To convert (\ref{eq:presentationofWperp_1}) into the target expression, where both of the original and the target expressions are regarded as those in $C$, the following three kinds of transformations are enough and an insertion of any subword $s_\beta s_\beta$ is \emph{not} required:
\begin{description}
\item[(T1)] a cancellation of a subword $s_\beta s_\beta$, with $s_\beta$ a loop generator;
\item[(T2)] a GBM that contains a loop generator;
\item[(T3)] a transformation $u \leadsto u'$, where $u$ and $u'$ contain no loop generators and $u \equiv u'$ in $Y$.
\end{description}
Now we proceed the proof by induction on the total number $N$ of transformations T1 and T2 in the conversion process.
Note that, in the case $N=0$, the original and the target expressions with respect to the generators of $W^{\perp [x]}$ are already equal, hence we have nothing to do.

In the conversion process of the general case, we apply (possibly no) T3 transformations to the expression (\ref{eq:presentationofWperp_1}) before applying the first T1 or T2 transformation, obtaining an expression
\begin{equation}
\label{eq:presentationofWperp_2}
u_1l_1u_2l_2u_3 \cdots u_nl_nu_{n+1}, \mbox{ where the } u_i \mbox{ contain no loop generators and } u_i \equiv w_{i-1}{}^{-1}w_i \enspace.
\end{equation}
First we consider the case that T1 is then applied to (\ref{eq:presentationofWperp_2}).
If the generators $l_{i-1}$ and $l_i$ with $2 \leq i \leq n$ are cancelled by the T1, then $u_i$ should be empty (hence $w_{i-1} \equiv w_i$) and $l_{i-1}=l_i$.
Now the T1 results in the expression
\begin{equation}
\label{eq:presentationofWperp_3}
u_1l_1u_2 \cdots u_{i-2}l_{i-2}u_{i-1}u_{i+1}l_{i+1}u_{i+2} \cdots u_nl_nu_{n+1} \enspace.
\end{equation}
Moreover, since now $u_{i-1}u_{i+1} \equiv {w_{i-2}}^{-1}w_{i+1}$, (\ref{eq:presentationofWperp_3}) is convertible by a sequence of T3 transformations into
\begin{equation}
\label{eq:presentationofWperp_4}
w_1l_1{w_1}^{-1} \cdots w_{i-2}l_{i-2}{w_{i-2}}^{-1}w_{i+1}l_{i+1}{w_{i+1}}^{-1} \cdots w_nl_n{w_n}^{-1} \enspace,
\end{equation}
which is also obtained from (\ref{eq:presentationofWperp_1}) by the transformation $r_x(w_i,l_i)^2 \leadsto 1$ in the statement (recall that now $w_{i-1} \equiv w_i$ in $Y$ and $l_{i-1}=l_i$, hence $r_x(w_{i-1},l_{i-1}) = r_x(w_i,l_i)$).
Now the total number of T1 and T2 transformations required in a conversion of (\ref{eq:presentationofWperp_4}) into the target expression is the same as that for (\ref{eq:presentationofWperp_3}), which is $N-1$, hence the claim follows by the induction assumption.

From now, we consider the other case that T2 is applied to (\ref{eq:presentationofWperp_2}) instead of T1.
Then the subword of (\ref{eq:presentationofWperp_2}) on which the GBM acts can be written as $bl_iu_{i+1} \cdots u_jl_jc$ for some $1 \leq i < j \leq n$, where we decomposed $u_i$ and $u_{j+1}$ as $u_i=ab$ and $u_{j+1}=cd$.
Here we only consider the case that the loop number $k=j-i+1$ of this GBM is odd, since the other case is similar.
Then Lemma \ref{lem:shuttlingandGBM} implies that, by the shape of shuttling tours (see Figure \ref{fig:shapeofcomponent}), there is a shuttling tour $p=lqLq^{-1}$ of order $k$, where $l$ and $L$ are loop generators, such that
\begin{eqnarray*}
l_i=l_{i+2}= \cdots =l_j=l \,,\, l_{i+1}=l_{i+3}= \cdots =l_{j-1}=L \enspace,\\
u_{i+1}=u_{i+3}= \cdots =u_{j-1}=q \,,\, u_{i+2}=u_{i+4}= \cdots =u_j=q^{-1}
\end{eqnarray*}
and the expression $q$ admits two decompositions of the form $q = cc'$ and $q = b^{-1}b'{}^{-1}$.
By these and the relations $u_m \equiv {w_{m-1}}^{-1}w_m$, we have $w_i \equiv w_{i+2} \equiv \cdots \equiv w_j$ and $w_{i+1} \equiv w_{i+3} \equiv \cdots \equiv w_{j-1}$ in $Y$.
Now this GBM, whose loop number is $k$, is of the form
\begin{displaymath}
bl_iu_{i+1} \cdots u_jl_jc=b\bigl(lqLq^{-1}\bigr)^{(k-1)/2}lc \overset{GBM}{\leadsto} b'{}^{-1}Lq^{-1}\bigl(lqLq^{-1}\bigr)^{(k-3)/2}lqLc'{}^{-1} \enspace,
\end{displaymath}
which converts (\ref{eq:presentationofWperp_2}) into
\begin{equation}
\label{eq:presentationofWperp_5}
u_1l_1 \cdots u_{i-1}l_{i-1}a\Bigl(b'{}^{-1}Lq^{-1}\bigl(lqLq^{-1}\bigr)^{(k-3)/2}lqLc'{}^{-1}\Bigr)dl_{j+1}u_{j+2} \cdots l_nu_{n+1} \enspace.
\end{equation}
On the other hand, the shuttling tour $p$ defines a relation $(w_i,l) \overset{k}{\sim} (w_iq,L)$, while $w_iq=w_iu_{i+1} \equiv w_{i+1}$.
Thus we have
\begin{eqnarray}
\nonumber
\mbox{(\ref{eq:presentationofWperp_1})} &\equiv&w_1l_1w_1^{-1} \cdots w_{i-1}l_{i-1}w_{i-1}{}^{-1}\bigl(w_ilw_i^{-1}w_{i+1}Lw_{i+1}{}^{-1}\bigr)^{(k-1)/2}w_ilw_i^{-1}\\
\nonumber
&&\quad \cdot\, w_{j+1}l_{j+1}w_{j+1}{}^{-1} \cdots w_nl_nw_n^{-1}\\
\nonumber
&\leadsto&w_1l_1w_1^{-1} \cdots w_{i-1}l_{i-1}w_{i-1}{}^{-1}\bigl(w_{i+1}Lw_{i+1}{}^{-1}w_ilw_i^{-1}\bigr)^{(k-1)/2}w_{i+1}Lw_{i+1}{}^{-1}\\
\label{eq:presentationofWperp_6}
&&\quad \cdot\, w_{j+1}l_{j+1}w_{j+1}{}^{-1} \cdots w_nl_nw_n^{-1} \enspace,
\end{eqnarray}
where the left-hand and the middle terms coincide with each other as expressions with respect to the generators of $W^{\perp [x]}$ given in the statement, and the last transformation comes from a transformation of second type given in the statement concerning the generators $r_x(w_i,l)$ and $r_x(w_{i+1},L)$.
Moreover, we have
\begin{displaymath}
w_{i-1}{}^{-1}w_{i+1} \equiv u_iu_{i+1} \equiv ab'{}^{-1},\ w_{i+1}{}^{-1}w_i \equiv q^{-1},\ w_{i+1}{}^{-1}w_{j+1} \equiv u_ju_{j+1} \equiv c'{}^{-1}d \enspace,
\end{displaymath}
therefore the right-hand side of (\ref{eq:presentationofWperp_6}) is convertible into (\ref{eq:presentationofWperp_5}) by a sequence of T3 transformations.
Thus the total number of T1 and T2 transformations required in a conversion of the right-hand side of (\ref{eq:presentationofWperp_6}) into the target expression is the same as that for (\ref{eq:presentationofWperp_5}), which is $N-1$, therefore the claim follows from the induction assumption.
Hence the proof of Theorem \ref{thm:presentationofWperp} is concluded.
\end{proof}
This theorem yields the main result of this subsection, as follows:
\begin{thm}
\label{thm:relationinWperp}
We have $r_x(w,\gamma)=r_x(u,\beta)$ if and only if $(w,\gamma) \sim (u,\beta)$.
Moreover, $r_x(w,\gamma)r_x(u,\beta)$ has order $k$ with $2 \leq k<\infty$ if and only if there exist $(w',\gamma')$ and $(u',\beta')$ in $\mathcal{R}_x$ such that $(w,\gamma) \sim (w',\gamma') \overset{k}{\sim} (u',\beta') \sim (u,\beta)$.
\end{thm}
\begin{proof}
Note that both of the two \lq\lq if'' parts follow from the definition.
From now, we prove the two \lq\lq only if'' parts.
First, note that if $(w,\gamma) \sim (w',\gamma')$, $(u,\beta) \sim (u',\beta')$, $(w,\gamma) \overset{k}{\sim} (u,\beta)$ and $(w',\gamma') \overset{\ell}{\sim} (u',\beta')$, then we have $k=\ell$, since now $k$ and $\ell$ are the order of the same element $r_x(w,\gamma)r_x(u,\beta)=r_x(w',\gamma')r_x(u',\beta')$.
Now Theorem \ref{thm:presentationofWperp} yields a presentation of $W^{\perp [x]}$ with generators given by
\begin{displaymath}
r_x\overline{(w,\gamma)}=r_x(w,\gamma)\,, \mbox{ where } \overline{(w,\gamma)} \mbox{ denotes an equivalence class of } \sim
\end{displaymath}
and fundamental relations given by
\begin{itemize}
\item $r_x\overline{(w,\gamma)}{}^2=1$ for every generator $r_x\overline{(w,\gamma)}$,
\item $\bigl(r_x\overline{(w,\gamma)}r_x\overline{(u,\beta)}\bigr)^k=1$ for every pair $(w,\gamma) \overset{k}{\sim} (u,\beta)$.
\end{itemize}
The above remark shows that these data form a Coxeter presentation of the Coxeter group $W^{\perp [x]}$.
Thus for two $(w,\gamma)$ and $(u,\beta)$ in $\mathcal{R}_x$, the product $r_x(w,\gamma)r_x(u,\beta)$ has order $1 \leq k<\infty$ \emph{if and only if} the relation \lq\lq $\bigl(r_x\overline{(w,\gamma)}r_x\overline{(u,\beta)}\bigr)^k=1$'' appears as one of the above fundamental relations.
Hence the claim follows.
\end{proof}
To study the relation $\sim$ further, we introduce an auxiliary graph $\mathcal{I}$ as follows.
The vertices of $\mathcal{I}$ are the loop generators $s_\gamma$ of $C$, or equivalently the loops in the graph $\mathcal{C}$.
The edges of $\mathcal{I}$ from a vertex $s_\beta$ to a vertex $s_\gamma$ are the paths $q$ in $\mathcal{Y}^1$ such that $s_\gamma qs_\beta q^{-1}$ is a shuttling tour of order one (note that such a path $q$ must be nonempty).
Let $\iota:\pi_1(\mathcal{I};*,*) \to \pi_1(\mathcal{Y}^1;*,*)$ be the groupoid homomorphism that maps a vertex $s_\gamma$ of $\mathcal{I}$ to $x \in \mathcal{V}(C)$ with $s_\gamma \in C_{x,x}$, and maps an edge $q$ of $\mathcal{I}$ to a path $q$ in $\mathcal{Y}^1$.
\begin{rem}
\label{rem:iotaisinjective}
Note that the induced homomorphism $\iota:\pi_1(\mathcal{I};s_\gamma) \to \pi_1(\mathcal{Y}^1;x)$ between fundamental groups, where $s_\gamma \in C_{x,x}$, is \emph{injective}.
Indeed, if a cancellation of two consecutive edges $ee^{-1}$ occurs in a path $\iota(q_n \cdots q_2q_1)$, then it should occur between the end of some $\iota(q_i)$ and the beginning of $\iota(q_{i+1})$, forcing $q_{i+1}=q_i^{-1}$ by definition of $\mathcal{I}$ and reducing $q_n \cdots q_1$ to a shorter path $q_n \cdots q_{i+2} q_{i-1} \cdots q_1$ in $\mathcal{I}$.
This argument implies the injectivity of $\iota$.
\end{rem}
Now Theorem \ref{thm:relationinWperp} implies that $r_x(w,\gamma)=r_x(u,\beta)$ if and only if there is a path $\widetilde{q}$ in $\mathcal{I}$ from $s_\beta$ to $s_\gamma$ such that $u=w\iota(\widetilde{q})_{(x)}$.
Hence $r_x(w,\gamma)=r_x(u,\gamma)$ if and only if $(w^{-1}u)_{(y)} \in \iota(\pi_1(\mathcal{I};s_\gamma))$ where $s_\gamma \in C_{y,y}$.

By this formulation, we obtain the following observation about the \emph{finite part} of the Coxeter group $W^{\perp I}$, denoted by $W^{\perp I}{}_{\mathrm{fin}}$, which we define as the product of all the irreducible components of $W^{\perp I}$ of finite type.
We prepare some temporary notations.
For a fixed loop generator $s_\gamma \in C_{x,x}$, put $P=\pi_1(\mathcal{I};s_\gamma)$ and let $N$ denote the kernel of the natural projection $\pi_1(\mathcal{Y}^1;x) \twoheadrightarrow \pi_1(\mathcal{Y};x)$ where the image of $c$ is denoted by $\overline{c}$.
Let $\overline{\iota}$ denote the composition $\pi_1(\mathcal{I};*,*) \overset{\iota}{\to} \pi_1(\mathcal{Y}^1;*,*) \twoheadrightarrow \pi_1(\mathcal{Y};*,*)$ of the above map $\iota$ followed by the natural projection.
Let $\mathrm{Ab}$ denote the abelianization map $\pi_1(\mathcal{Y}^1;x) \twoheadrightarrow \mathrm{Ab}(\pi_1(\mathcal{Y}^1;x))$.
\begin{thm}
\label{thm:conditionfornonfinitepart}
Under the above notations, suppose that $|S|<\infty$ and the subgroup $\mathrm{Ab}(\iota(P)N)$ has infinite index in $\mathrm{Ab}(\pi_1(\mathcal{Y}^1;x))$.
Then for any $w \in Y_I$, the vertex $r(w,\gamma)$ of the Coxeter graph of $W^{\perp I}$ is adjacent to an edge with label $\infty$, hence $r(w,\gamma) \not\in W^{\perp I}{}_{\mathrm{fin}}$.

In particular, the second assumption is satisfied if $\mathrm{rk}(P)+n<\mathrm{rk}(\pi_1(\mathcal{Y}^1;x))$, where $\mathrm{rk}(F)$ denotes the rank of a free group $F$ and $n$ is the number of $2$-cells of $\mathcal{Y}$.
\end{thm}
\begin{proof}
To prove the first claim, we suppose that $r(w,\gamma)$ is not adjacent to an edge with label $\infty$ in the Coxeter graph of $W^{\perp I}$, and give a decomposition of $\mathrm{Ab}(\pi_1(\mathcal{Y}^1;x))$ into a finite number of cosets modulo $\mathrm{Ab}(\iota(P)N)$, yielding a contradiction.
Theorem \ref{thm:decompofC}(4) allows us to assume without loss of generality that $w=1$.
Then for any $a \in \pi_1(\mathcal{Y}^1;x)$, we have either $r(1,\gamma)=r(\overline{a}_{(x_I)},\gamma)$, or $r(1,\gamma)r(\overline{a}_{(x_I)},\gamma)$ has order $2 \leq k<\infty$, by the above hypothesis of the proof.

In the second case, Theorem \ref{thm:relationinWperp} implies the existence of an element $b \in \pi_1(\mathcal{Y}^1;x_I)$ and a shuttling tour $s_\beta qs_{\beta'}q^{-1}$ of order $k$ such that
\begin{equation}
\label{eq:conditionfornonfinitepart_1}
(1,\gamma) \sim (\overline{b},\beta) \overset{k}{\sim} (\overline{b}(\overline{q})_{(x_I)},\beta') \sim (\overline{a}_{(x_I)},\gamma) \enspace.
\end{equation}
Now the hypothesis $|S|<\infty$ implies that only a finite number of shuttling tours, say $s_{\beta_i}q_is_{\beta'_i}q_i^{-1}$ with $1 \leq i \leq r<\infty$, appear in this manner.
Let the above-mentioned shuttling tour be the $j$-th one.
Then (\ref{eq:conditionfornonfinitepart_1}) implies that $\overline{b}=\overline{\iota}(\widetilde{q})_{(x_I)}$ for some $\widetilde{q} \in \pi_1(\mathcal{I};s_\gamma,s_{\beta_j})$, and $\overline{a}_{(x_I)}=\overline{b}(\overline{q_j})_{(x_I)}\overline{\iota}(\widetilde{q'})_{(x_I)}$ for some $\widetilde{q'} \in \pi_1(\mathcal{I};s_{\beta'_j},s_\gamma)$.
Thus we have $\overline{a}_{(x_I)}=\overline{(\iota(\widetilde{q})q_j\iota(\widetilde{q'}))_{(x_I)}}$.
Since $s_{\gamma} \in C_{x,x}$ and $a \in \pi_1(\mathcal{Y}^1;x)$, it follows that
\begin{displaymath}
\overline{a}=\overline{\iota(\widetilde{q})q_j\iota(\widetilde{q'})} \in \overline{\iota(P)\iota(\widetilde{p_j})q_j\iota(\widetilde{p'_j})\iota(P)}
\end{displaymath}
for some $\widetilde{p_j} \in \pi_1(\mathcal{I};s_\gamma,s_{\beta_j})$ and $\widetilde{p'_j} \in \pi_1(\mathcal{I};s_{\beta'_j},s_\gamma)$ (fixed for each $j$), therefore $a \in \iota(P)c_j\iota(P)N$ where we put $c_j=\iota(\widetilde{p_j})q_j\iota(\widetilde{p'_j}) \in \pi_1(\mathcal{Y}^1;x)$.

On the other hand, in the first case $r(1,\gamma)=r(\overline{a}_{(x_I)},\gamma)$, we have $\overline{a} \in \overline{\iota}(P)$ and $a \in \iota(P)c_0\iota(P)N$, where we put $c_0=1 \in \pi_1(\mathcal{Y}^1;x)$.

Summarizing, we have $a \in \iota(P)c_i\iota(P)N$ for some $0 \leq i \leq r$.
Since $a$ was arbitrarily chosen, this means that $\pi_1(\mathcal{Y}^1;x)=\bigcup_{i=0}^{r}\iota(P)c_i\iota(P)N$, therefore we have
\begin{displaymath}
\mathrm{Ab}(\pi_1(\mathcal{Y}^1;x))=\bigcup_{i=0}^{r}\mathrm{Ab}(\iota(P))\mathrm{Ab}(c_i)\mathrm{Ab}(\iota(P))\mathrm{Ab}(N)=\bigcup_{i=0}^{r}\mathrm{Ab}(c_i)\mathrm{Ab}(\iota(P)N) \enspace.
\end{displaymath}
This is the desired decomposition of $\mathrm{Ab}(\pi_1(\mathcal{Y}^1;x))$ into finitely many cosets modulo $\mathrm{Ab}(\iota(P)N)$.
Hence the first claim holds.

For the second claim, let $\widetilde{g}_1,\dots,\widetilde{g}_r$ (where $r=\mathrm{rk}(P)$) be free generators of $P$ and $c_1,\dots,c_n$ the boundaries of the $2$-cells of $\mathcal{Y}$.
Then $N$ is generated by the conjugates of the $(c_i)_{(x)}$, therefore $\mathrm{Ab}(\iota(P)N)$ is generated by the $\mathrm{Ab}(\iota(\widetilde{g}_i))$ and $\mathrm{Ab}((c_j)_{(x)})$.
By the hypothesis, the total number $r+n$ of these generators is less than $\mathrm{rk}(\pi_1(\mathcal{Y}^1;x))$ that is equal to the rank of the free abelian group $\mathrm{Ab}(\pi_1(\mathcal{Y}^1;x))$, implying that the subgroup $\mathrm{Ab}(\iota(P)N) \subseteq \mathrm{Ab}(\pi_1(\mathcal{Y}^1;x))$ has infinite index, as desired.
Hence Theorem \ref{thm:conditionfornonfinitepart} holds.
\end{proof}
On the other hand, for a general $S$, we have the following result.
For any $J \subseteq S$, put
\begin{equation}
\label{eq:defofJperp}
J^\perp=\{s \in S \smallsetminus J \mid s \mbox{ is not adjacent to } J \mbox{ in the graph } \Gamma\} \enspace.
\end{equation}
\begin{prop}
\label{prop:relationbetweenfiniteparts}
Let $s_\gamma=w_x^s$ be a loop generator of $C$, $w \in Y_I$ and $J \subseteq S$ a finite subset such that $[x]_{\sim s} \subseteq J$ and $[x] \smallsetminus J \subseteq J^\perp$, hence $s_\gamma \in R^{J,[x] \cap J}$.
Then we have $r(w,\gamma) \not\in W^{\perp I}{}_{\mathrm{fin}}$ whenever $s_\gamma \not\in W_J^{\perp [x] \cap J}{}_{\mathrm{fin}}$.
\end{prop}
\begin{proof}
First, Theorem \ref{thm:decompofC}(4) allows us to assume that $w=1$ and $x=x_I$, therefore $r(w,\gamma)=s_\gamma$.
Now since $[x] \smallsetminus J \subseteq J^\perp$, we have $\Phi_J^{\perp [x]}=\Phi_J^{\perp [x] \cap J}$, hence $\Pi^{J,[x] \cap J}=\Pi^{J,[x]} \subseteq \Pi^{[x]}$.
This implies that the Coxeter system $(W_J^{\perp [x] \cap J},R^{J,[x] \cap J})$ is a subsystem of $(W^{\perp [x]},R^{[x]})$.
Hence if $s_\gamma \not\in W_J^{\perp [x] \cap J}{}_{\mathrm{fin}}$, then the irreducible component of $W^{\perp [x]}$ containing $s_\gamma=r(w,\gamma)$ contains the infinite irreducible component of $W_J^{\perp [x] \cap J}$ containing $s_\gamma$, implying that $r(w,\gamma) \not\in W^{\perp I}{}_{\mathrm{fin}}$ as desired.
\end{proof}
At the last of this subsection, we explain how to compute the order of a shuttling tour in an individual case.
We fix a shuttling tour $s_\gamma us_\beta u^{-1}$ in a connected component $\mathcal{G}$ of $\mathcal{C}(J)$ (recall the definition in Section \ref{sec:graphC}) with $\beta,\gamma \in \Phi^+$, and put $s_\beta=w_x^s$, $s_\gamma=w_y^t$ and $J=[x] \cup \{s,s'\}=[y] \cup \{t,t'\}$.
Put $K=J_{\sim (J \smallsetminus [y])}$, which is of finite type.
Now by applying our arguments above to $W_J^{\perp [y]}$ instead of $W^{\perp [y]}$, it follows that $R^{J,[y]}=\{s_\gamma,us_\beta u^{-1}=s_{u \cdot \beta}\}$ and, since $u \cdot \beta \in \Phi^+$, we have $\Pi^{J,[y]}=\{\gamma, u \cdot \beta\}$.
Thus Theorem \ref{thm:conditionforrootbasis} says that the order $k$ of the shuttling tour $s_\gamma us_\beta u^{-1}$ is determined by the equality
\begin{displaymath}
\langle \gamma, u \cdot \beta \rangle=-\cos(\pi/k) \enspace.
\end{displaymath}
Moreover, we have $\Phi_J^{\perp [y]}=\Phi_K^{\perp [y] \cap K}$ by definition of $K$, while $W_J^{\perp [y]}$ is now the dihedral group of order $2k$, having $k$ positive roots.
Thus we have another relation
\begin{displaymath}
|(\Phi_J^{\perp [y]})^+|=|(\Phi_K^{\perp [y] \cap K})^+|=k \enspace,
\end{displaymath}
which also enables us to determine the order $k$ if the root system $\Phi_K$ of a finite Coxeter group $W_K$ is well understood.

Owing to these arguments, if $W_K$ is not irreducible, namely if $t$ and $t'$ lie in distinct irreducible components of $W_J$, then we have
\begin{displaymath}
k=1 \mbox{ if } w_y^{t'} \cdot y \neq y \,, \mbox{ and } k=2 \mbox{ if } w_y^{t'} \cdot y=y \enspace.
\end{displaymath}
If $W_K$ is irreducible, then the possible situations are completely listed (up to symmetry) in Tables \ref{tab:listofedges_1} and \ref{tab:listofedges_2}, where we use the canonical labelling $K=\{r_1,\dots,r_n\}$ given in Section \ref{sec:longestelement}, and abbreviate $r_i$ to $i$.
%%%%%
\begin{table}[htb]
\centering
\caption{List of the shuttling tours, in the case $\mathrm{Type}\,K=A_n,B_n$ or $D_n$}
\label{tab:listofedges_1}
\begin{tabular}{cc|cc|c}
$K$ && $(s,s')$ & $(t,t')$ & order\\ \hline
$A_n$ & ($n \geq 3$) & $(1,2)$ & $(n,n-1)$ & $1$\\
& ($n=2$) & $(1,2)$ & $(2,1)$ & $3$\\ \hline
$B_n$ && $(1,2)$ & $(2,1)$ & $4$\\
&& $(3,1)$ & $(3,2)$ & $2$\\
&& $(4,2)$ & $(4,2)$ & $2$\\
& ($i \geq 5$) & $(i,i-2)$ & $(i,i-2)$ & $1$\\
& ($i \geq 4$) & $(i,i-1)$ & $(i,i-1)$ & $1$\\ \hline
$D_n$ && $(1,2)$ & $(1,2)$ & $2$\\
& ($n \geq 6$) & $(4,2)$ & $(4,2)$ & $2$\\
& ($4 \neq i \leq n-2$) & $(i,i-2)$ & $(i,i-2)$ & $1$\\
& ($n \geq 5$ even) & $(n-1,n-2)$ & $(n-1,n-2)$ & $1$\\
& ($n \geq 5$ odd) & $(n-1,n-2)$ & $(n,n-2)$ & $1$
\end{tabular}
\end{table}
%%%%%
\begin{table}[htb]
\centering
\caption{List of the shuttling tours, for the remaining cases}
\label{tab:listofedges_2}
\begin{tabular}{c|cc|c||c|cc|c}
$K$ & $(s,s')$ & $(t,t')$ & order & $K$ & $(s,s')$ & $(t,t')$ & order\\ \hline
$E_6$ & $(1,3)$ & $(6,5)$ & $1$ & $E_8$ & $(1,3)$ & $(1,3)$ & $1$\\
& $(2,4)$ & $(2,4)$ & $3$ && $(1,8)$ & $(1,8)$ & $2$\\
& $(3,6)$ & $(5,1)$ & $1$ && $(2,4)$ & $(2,4)$ & $1$\\ \cline{1-4}
$E_7$ & $(1,3)$ & $(1,3)$ & $3$ && $(2,7)$ & $(2,7)$ & $1$\\
& $(2,4)$ & $(2,4)$ & $1$ && $(3,6)$ & $(3,6)$ & $1$\\
& $(2,7)$ & $(2,7)$ & $1$ && $(7,1)$ & $(7,1)$ & $1$\\
& $(3,6)$ & $(3,6)$ & $1$ && $(8,7)$ & $(8,7)$ & $3$\\ \cline{5-8}
& $(6,1)$ & $(6,1)$ & $2$ & $H_3$ & $(1,2)$ & $(3,2)$ & $2$\\ \cline{5-8}
& $(7,6)$ & $(7,6)$ & $1$ & $H_4$ & $(1,2)$ & $(1,2)$ & $3$\\ \cline{1-4}
$F_4$ & $(1,2)$ & $(1,2)$ & $3$ && $(2,4)$ & $(2,4)$ & $2$\\
& $(1,4)$ & $(4,1)$ & $4$ && $(4,3)$ & $(4,3)$ & $5$\\ \cline{5-8}
& $(2,4)$ & $(3,1)$ & $2$ & $I_2(m)$ & $(1,2)$ & $(2,1)$ & $m$\\
\end{tabular}
\end{table}
%%%%%
%%%%%
\section{Description of the entire centralizer}
\label{sec:descriptionoffullZ}
Although the subgroup $C_I=C_{x_I,x_I}$ we have described occupies a fairly large part of the centralizer $Z_W(W_I)$, it is in fact not yet the whole of $Z_W(W_I)$.
In this section, we describe the structure of the entire centralizer $Z_W(W_I)$ further (Section \ref{sec:descriptionoffullZ_centralizer}).
On the other hand, by using the above results, we can also decompose the normalizer $N_W(W_I)$ in a different manner from \cite{Bri-How}.
We also show the new decomposition of $N_W(W_I)$ as a by-product of our argument (Section \ref{sec:descriptionoffullZ_normalizer}).
It is worthy to notice that our decomposition of $N_W(W_I)$ possesses some good properties which failed in \cite{Bri-How}; for instance, in contrast with the second factor $\widetilde{N}_I$ of the decomposition of $N_W(W_I)$ in \cite{Bri-How} which is a Coxeter group as well, our second factor $W^{\perp I}$ is a reflection subgroup of $W$, not just a general Coxeter group, and the factor admits a simple closed definition.
%%%%%
\subsection{On the centralizers}
\label{sec:descriptionoffullZ_centralizer}
The aim here is to show that $Z_W(W_I)$ is a (not necessarily split) extension of $C_I$ by an elementary abelian $2$-group $\widetilde{\mathcal{A}}$ described below, and give a further decomposition of $Z_W(W_I)$.
Our first observation is that the map $x_\lambda \mapsto y_\lambda$ defines an isomorphism $W_{[x]} \to W_{[y]}$ for any $x,y \in \mathcal{V}(C)$, since this map is induced by conjugation by an element of $C_{y,x} \neq \emptyset$.
Hence if $A \subseteq \Lambda$, the groups $W_{x_A}$ are isomorphic for all $x \in \mathcal{V}(C)$ in this manner, justifying the expression \lq\lq $A$ is an irreducible component of $\Lambda$'' to mean that $W_{x_A}$ is an irreducible component of $W_{[x]}$.
The situation is similar for the expressions \lq\lq $A \subseteq \Lambda$ is of finite type'' and \lq\lq $A \subseteq \Lambda$ is of $(-1)$-type''.
Note that, if $A \subseteq \Lambda$ is of finite type, then
\begin{displaymath}
ww_0(x_A)w^{-1}=w_0(y_A) \mbox{ holds in } W \mbox{ for all } w \in C_{y,x}
\end{displaymath}
since the above isomorphism $W_{x_A} \to W_{y_A}$ sends $w_0(x_A)$ to $w_0(y_A)$.

Let $A \subseteq \Lambda$ be of finite type and a union of irreducible components of $\Lambda$.
Then define
\begin{displaymath}
x^A \in S^{(\Lambda)} \mbox{ for } x \in \mathcal{V}(C) \mbox{ by } x^A{}_\lambda=w_0(x_A)x_\lambda w_0(x_A) \enspace.
\end{displaymath}
Note that there is a common permutation $\sigma_A$ on $\Lambda$ such that $x^A{}_\lambda=x_{\sigma_A(\lambda)}$ for all $x \in \mathcal{V}(C)$ and $\lambda \in \Lambda$, and that $w_0(x_A)$ and $w_0(x_{A'})$ commute (hence $\sigma_A$ and $\sigma_{A'}$ do as well) for such $A,A' \subseteq \Lambda$.
More precisely, $w_0(x_A)w_0(x_{A'})=w_0(x_{AA'})$ where $AA'$ denotes the symmetric difference of $A$ and $A'$ (hence $AA'=A'A$).
Let
\begin{displaymath}
\widetilde{\mathcal{A}}=\{A \subseteq \Lambda \mid x_I{}^A \mbox{ is defined and } x_I{}^A \in \mathcal{V}(C)\} \enspace,
\end{displaymath}
which is an elementary abelian $2$-group with the symmetric difference as multiplication.
Indeed, if $A,A' \in \widetilde{\mathcal{A}}$ and $w \in C_{x_I{}^A,x_I} \neq \emptyset$, then $w$ also lies in $C_{(x_I{}^A)^{A'},x_I{}^{A'}}$ and $x_I{}^{A'} \in \mathcal{V}(C)$, hence $(x_I{}^A)^{A'}=x_I{}^{AA'} \in \mathcal{V}(C)$ and $AA' \in \widetilde{\mathcal{A}}$.
More precisely, for $A \in \widetilde{\mathcal{A}}$, the graph $\mathcal{C}$ defined in Section \ref{sec:graphC} admits an automorphism $\tau_A$ such that
\begin{eqnarray}
\nonumber
&&\tau_A(y)=y^A \mbox{ for } y \in \mathcal{V}(C)\,,\, \tau_A(w_y^s)=w_{y^A}^s \mbox{ for any edge } w_y^s \enspace,\\
\label{eq:tauA}
&&w_0(z_A)ww_0(y_A)=\tau_A(w) \mbox{ holds in } C_{z^A,y^A} \mbox{ for all } w \in C_{z,y} \enspace.
\end{eqnarray}
Note that the map $A \mapsto \tau_A$ is a group homomorphism $\widetilde{\mathcal{A}} \to \mathrm{Aut}\,\mathcal{C}$.

\begin{figure}[htb]
\centering
\begin{picture}(140,90)
\put(10,35){\circle*{10}}\put(10,20){\hbox to0pt{\hss$s_1$\hss}}
\put(50,35){\circle*{10}}\put(50,20){\hbox to0pt{\hss$s_2$\hss}}
\put(90,15){\circle{10}}\put(90,0){\hbox to0pt{\hss$s_3$\hss}}
\put(130,15){\circle*{10}}\put(130,0){\hbox to0pt{\hss$s_4$\hss}}
\put(130,55){\circle*{10}}\put(130,65){\hbox to0pt{\hss$s_5$\hss}}
\put(90,55){\circle{10}}\put(90,65){\hbox to0pt{\hss$s_6$\hss}}
\put(15,35){\line(1,0){30}}
\put(55,37){\line(2,1){30}}
\put(55,33){\line(2,-1){30}}
\put(95,15){\line(1,0){30}}
\put(95,55){\line(1,0){30}}
\put(130,20){\line(0,1){30}}
\put(70,80){\hbox to0pt{\hss (A)\hss}}
\end{picture}
\qquad\qquad
\begin{picture}(120,90)
\put(20,50){\hbox to0pt{\hss$(1,2,4,5)$\hss}}
\put(20,10){\hbox to0pt{\hss$(5,4,2,1)$\hss}}
\put(100,10){\hbox to0pt{\hss$(2,1,5,4)$\hss}}
\put(100,50){\hbox to0pt{\hss$(4,5,1,2)$\hss}}
\put(20,20){\line(0,1){25}}\put(13,30){\hbox to0pt{\hss$6$\hss}}
\put(100,20){\line(0,1){25}}\put(108,30){\hbox to0pt{\hss$6$\hss}}
\put(45,14){\line(1,0){30}}\put(60,2){\hbox to0pt{\hss$3$\hss}}
\multiput(45,54)(12,0){3}{\line(1,0){6}}\put(60,60){\hbox to0pt{\hss$3$\hss}}
\put(60,80){\hbox to0pt{\hss (B)\hss}}
\end{picture}
\caption{(A) Coxeter graph for Example \ref{ex:B_not_split} (here black circles signify elements of $I$); (B) the resulting $1$-skeleton $\mathcal{Y}^1$, with maximal tree $\mathcal{T}$ depicted by solid lines}
\label{fig:Coxetergraph_B_not_split}
\end{figure}
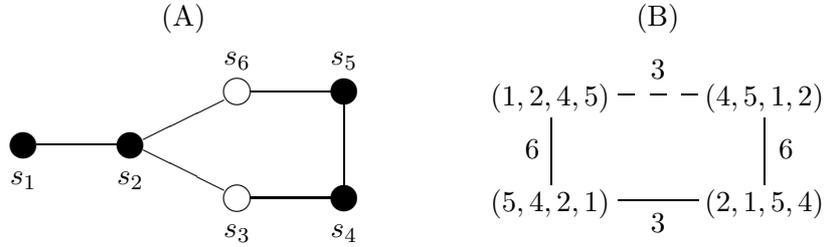

Now we have the following result.
Recall the maximal tree $\mathcal{T}$ in $\mathcal{Y}^1$ and the paths $p_{y,x}$ in $\mathcal{T}$ introduced in Section \ref{sec:decompofC_factorY}.
Note that each path in $\mathcal{C}$ from $x$ to $y$ represents an element of $C_{y,x}$.
\begin{thm}
\label{thm:exactsequence_Z}
We have an exact sequence $1 \to C_I \hookrightarrow Z_W(W_I) \overset{\mathsf{A}}{\to} \widetilde{\mathcal{A}} \to 1$ of group homomorphisms and a map $g:\widetilde{\mathcal{A}} \to Z_W(W_I)$ with $\mathsf{A} \circ g=\mathrm{id}_{\widetilde{\mathcal{A}}}$, where
\begin{eqnarray*}
&&\mathsf{A}:Z_W(W_I) \to \widetilde{\mathcal{A}}\,,\, w \mapsto \mathsf{A}(w)=\{\lambda \in \Lambda \mid w \cdot \alpha_{(x_I)_\lambda}=-\alpha_{(x_I)_\lambda}\} \enspace,\\
&&g:\widetilde{\mathcal{A}} \to Z_W(W_I)\,,\, A \mapsto g_A=p_{x_I,x_I{}^A}w_0((x_I)_A) \enspace.
\end{eqnarray*}

If we define a subgroup $\widetilde{B}_I$ of $Z_W(W_I)$ by
\begin{displaymath}
\widetilde{B}_I=\{w \in Z_W(W_I) \mid (\Phi^{\perp I})^+=w \cdot (\Phi^{\perp I})^+\} \enspace,
\end{displaymath}
then we have $g(\widetilde{\mathcal{A}}) \subseteq \widetilde{B}_I$, the sequence $1 \to Y_I \hookrightarrow \widetilde{B}_I \overset{\mathsf{A}}{\to} \widetilde{\mathcal{A}} \to 1$ is exact and $Z_W(W_I)=W^{\perp I} \rtimes \widetilde{B}_I$, where $\widetilde{B}_I$ acts on $W^{\perp I}$ as automorphisms of Coxeter system.

Moreover, if $p \in \pi_1(\mathcal{Y}^1;x_I)$ and $A,A' \in \widetilde{\mathcal{A}}$, then we have
\begin{eqnarray}
\label{eq:relationIofG}
g_Apg_A{}^{-1}&=&\tau_A(p)_{(x_I)} \enspace,\\
\label{eq:relationIIofG}
g_Ag_{A'}&=&\tau_A(p_{x_I,x_I{}^{A'}})_{(x_I)}g_{AA'} \enspace,\\
\label{eq:relationIIIofG}
g_A{}^2&=&\tau_A(p_{x_I,x_I{}^A})_{(x_I)} \enspace,\\
\label{eq:relationIVofG}
g_Ag_{A'}g_A{}^{-1}g_{A'}{}^{-1}&=&\tau_A(p_{x_I,x_I{}^{A'}})_{(x_I)}\tau_{A'}(p_{x_I{}^A,x_I})_{(x_I)} \enspace.
\end{eqnarray}
\end{thm}
\begin{proof}
For the first part, the only nontrivial claim is the one on well-definedness of the homomorphism $\mathsf{A}$.
This will follow easily once we show that $\mathsf{A}(w) \in \widetilde{\mathcal{A}}$.
Now Lemma \ref{lem:rightdivisor} implies that $\mathsf{A}(w)$ is of finite type.
On the other hand, for $\lambda,\mu \in \Lambda$ with $\langle \alpha_{(x_I)_\lambda},\alpha_{(x_I)_\mu} \rangle \neq 0$, we have $\lambda \in \mathsf{A}(w)$ if and only if $\mu \in \mathsf{A}(w)$, since $w$ leaves the form $\langle \,,\,\rangle$ invariant.
Thus $\mathsf{A}(w)$ is a union of irreducible components of $\Lambda$, hence $x_I{}^{\mathsf{A}(w)}$ is defined.
Moreover, the remaining claim $x_I{}^{\mathsf{A}(w)} \in \mathcal{V}(C)$ follows since $w_0((x_I)_{\mathsf{A}(w)})w \in C_{x_I{}^{\mathsf{A}(w)},x_I}$.
Hence the first part is proven.

For the second part, we have $g_A \in \widetilde{B}_I$ for all $A \in \widetilde{\mathcal{A}}$, since $p_{x_I,x_I{}^A} \in Y$ and $w_0((x_I)_A) \in W_I$.
This implies the exactness of the sequence.
Moreover, the first exact sequence and the properties $g(\widetilde{\mathcal{A}}) \subseteq \widetilde{B}_I$ and $C_I=W^{\perp I}\,Y_I$ imply that $Z_W(W_I)=W^{\perp I}\widetilde{B}_I$, and then we have $Z_W(W_I)=W^{\perp I} \rtimes \widetilde{B}_I$ in the same way as the case of the decomposition $C_I=W^{\perp I} \rtimes Y_I$.
The conjugation action of $\widetilde{B}_I$ on $W^{\perp I}$ defines automorphisms of Coxeter system since each element of $\widetilde{B}_I$ leaves the positive system $(\Phi^{\perp I})^+$ of $W^{\perp I}$ invariant.

Finally, since $w_0(x_{A'})=w_0(x_A)w_0(x_{AA'})=w_0((x^{A'})_A)w_0(x_{AA'})$, the last formulae are deduced by straightforward computations based on (\ref{eq:tauA}).
Hence Theorem \ref{thm:exactsequence_Z} holds.
\end{proof}
For a further study, we define
\begin{eqnarray*}
\mathcal{A}&=&\{A \in \widetilde{\mathcal{A}} \mid A \cap A'=\emptyset \mbox{ if } A' \subseteq \Lambda \mbox{ is an irreducible component of } (-1)\mbox{-type}\} \enspace,\\
\mathcal{A}'&=&\{A \subseteq \Lambda \mid A \mbox{ is a finite union of irreducible components of } (-1)\mbox{-type}\} \enspace.
\end{eqnarray*}
Since $x_I{}^A=x_I$ and $g_A=w_0((x_I)_A)$ for any $A \in \mathcal{A}'$, both $\mathcal{A}$ and $\mathcal{A}'$ are subgroups of $\widetilde{\mathcal{A}}$, and we have $\widetilde{\mathcal{A}}=\mathcal{A} \times \mathcal{A}'$ and $g(\mathcal{A}')=Z(W_I)$ by the structure of $Z(W_I)$ (see Section \ref{sec:longestelement}).
Moreover, $g_{AA'}=g_Ag_{A'}$ for $A \in \mathcal{A}$ and $A' \in \mathcal{A}'$.
Now the following result is an easy consequence of these properties and Theorem \ref{thm:exactsequence_Z}:
\begin{thm}
\label{thm:exactsequence_B}
If we define a subgroup $B_I$ of $\widetilde{B}_I$ by
\begin{displaymath}
B_I=\{w \in \widetilde{B}_I \mid \mathsf{A}(w) \in \mathcal{A}\} \enspace,
\end{displaymath}
then we have $g(\mathcal{A}) \subseteq B_I$, the sequence $1 \to Y_I \hookrightarrow B_I \overset{\mathsf{A}}{\to} \mathcal{A} \to 1$ is exact and $\widetilde{B}_I=Z(W_I) \times B_I$.
Moreover, we have
\begin{displaymath}
Z_W(W_I)=W^{\perp I} \rtimes (Z(W_I) \times B_I)=(Z(W_I) \times W^{\perp I}) \rtimes B_I \enspace,
\end{displaymath}
where $B_I$ acts on $W^{\perp I}$ as automorphisms of Coxeter system.
\end{thm}
Note that the exact sequences in these theorems may \emph{not} split, since the map $g$ is \emph{not} necessarily a group homomorphism (see below for a counterexample).
However, we have the following result:
\begin{prop}
\label{prop:Zsplits}
If the maximal tree $\mathcal{T} \subseteq \mathcal{Y}^1$ is stable under the maps $\tau_A$ for all $A \in \mathcal{A}$, then $g:\widetilde{\mathcal{A}} \to Z_W(W_I)$ is a group homomorphism.
Hence in this case, the exact sequences of Theorems \ref{thm:exactsequence_Z} and \ref{thm:exactsequence_B} split via the $g$.
\end{prop}
\begin{proof}
First, by definition of $\mathcal{A}'$, the $\tau_A$ is identity if $A \in \mathcal{A}'$.
Then, since $\widetilde{\mathcal{A}}=\mathcal{A} \times \mathcal{A}'$, the hypothesis implies that $\mathcal{T}$ is stable under the $\tau_A$ for all $A \in \widetilde{\mathcal{A}}$.
Thus $\tau_A(p_{x_I,x_I{}^{A'}})=p_{x_I{}^A,x_I{}^{AA'}}$ for all $A,A' \in \widetilde{\mathcal{A}}$, hence $g_Ag_{A'}=g_{AA'}$ by (\ref{eq:relationIIofG}) as desired.
\end{proof}
\begin{exmp}
\label{ex:B_not_split}
Here we give an example to show that the conclusion of Proposition \ref{prop:Zsplits} does not necessarily hold in general.
Let $(W,S)$ be the Coxeter system corresponding to the Coxeter graph in the left (A) of Figure \ref{fig:Coxetergraph_B_not_split}, and put $I=\{s_1,s_2,s_4,s_5\}$ and $x_I=(s_1,s_2,s_4,s_5)$ (therefore $\Lambda = \{1,2,3,4\}$).
Then a direct calculation shows that the graph $\mathcal{Y}^1$ is as in the right (B) of Figure \ref{fig:Coxetergraph_B_not_split}, where each vertex $(s_{i_1},s_{i_2},s_{i_3},s_{i_4})$ of $\mathcal{Y}^1$ is abbreviated to $(i_1,i_2,i_3,i_4)$ and each edge $w_y^s$ from a vertex $y$ to another vertex $w_y^s \cdot y$ is labelled $s$.
(See Section \ref{sec:example} for a more detailed example of calculation of the graph $\mathcal{Y}^1$.)
We take a maximal tree $\mathcal{T}$ in $\mathcal{Y}^1$ to be obtained by removing the edge between $x_I = (1,2,4,5)$ and $(4,5,1,2)$.
Then we have $\mathcal{A} = \widetilde{\mathcal{A}} = \{\emptyset,\Lambda\}$.
By (\ref{eq:relationIIIofG}), the element $g = g_{\Lambda}$ of $B_I$ corresponding to the generator $\Lambda$ of $\mathcal{A}$ satisfies that $g^2$ is the cycle $(1,2,4,5) \leftarrow (5,4,2,1) \leftarrow (2,1,5,4) \leftarrow (4,5,1,2) \leftarrow (1,2,4,5)$ in $\mathcal{Y}^1$, which generates the free group $Y_I \simeq \pi_1(\mathcal{Y};x_I) \simeq \mathbb{Z}$ of rank one.
This implies that $B_I$ is now the infinite cyclic group generated by $g$, in which $Y_I$ is a normal subgroup of index two.
Hence this $B_I$ does not split over $Y_I$, since $B_I$ has no subgroup isomorphic to $B_I/Y_I \simeq \mathbb{Z}/2\mathbb{Z}$.
\end{exmp}
%%%%%
%%%%%
A general result in combinatorial group theory (see e.g., \cite[Proposition 10.1]{Joh}) enables us to construct a presentation of a group extension from those of its factors.
In the case of the group $B_I$, this result gives the following presentation:
\begin{thm}
\label{thm:presentationofB}
Let $\mathcal{A}_0$ be a basis of the elementary abelian $2$-group $\mathcal{A}$.
Then $B_I$ admits a presentation, with generators given by $(w_\xi)_{(x_I)}$ for all $w_\xi \in E(\mathcal{Y})$ and $g_A$ for all $A \in \mathcal{A}_0$, and fundamental relations given by those of $Y_I$, relations (\ref{eq:relationIofG}) for $A \in \mathcal{A}_0$ and $q=(w_{\xi})_{(x_I)}$ (where $w_\xi \in E(\mathcal{Y})$), relations (\ref{eq:relationIIIofG}) for $A \in \mathcal{A}_0$, and relations (\ref{eq:relationIVofG}) for $A,A' \in \mathcal{A}_0$.
\end{thm}
Finally, we determine the action of $B_I$ on $W^{\perp I}$ in the following manner:
\begin{prop}
\label{prop:actiononWperp_B}
Let $w_y^s=s_{\gamma(y,s)}$ be a loop generator of $C$, where the root $\gamma(x,s)$ is the one defined in Lemma \ref{lem:charofBphi}, and let $w,u \in Y_I$ and $A \in \mathcal{A}$.
Then
\begin{displaymath}
ur(w,w_y^s)u^{-1}=r(uw,w_y^s) \mbox{ and } g_Ar(w,w_y^s)g_A{}^{-1}=r(\tau_A(wp_{x_I,y})_{(x_I)},w_{y^A}^s) \enspace.
\end{displaymath}
\end{prop}
\begin{proof}
Since $\tau_A(s_{\gamma(y,s)})=s_{\gamma(y^A,s)}$, this claim is an easy consequence of (\ref{eq:tauA}).
\end{proof}
%%
%%%%%
\subsection{On the normalizers}
\label{sec:descriptionoffullZ_normalizer}
In a similar manner, we also give a decomposition of $N_W(W_I)$, in particular that of the second factor $N_I$ of $N_W(W_I)=W_I \rtimes N_I$ (see Section \ref{sec:precedingonN}).

Let $\mathcal{A}_\mathrm{N}$ be the set of the bijections $\rho:\Lambda \to \Lambda$ such that $\rho(x_I) \in \mathcal{V}(C)$, where $\rho$ acts on $S^{(\Lambda)}$ by $\rho(y)_\lambda=y_{\rho^{-1}(\lambda)}$.
Then we have $N_I=\bigsqcup_{\rho \in \mathcal{A}_\mathrm{N}}C_{x_I,\rho(x_I)}$ by definition.
Moreover, the connectedness of $\mathcal{Y}^1$ implies that each $\rho$ fixes all but finitely many elements of $\Lambda$, hence $\rho$ has finite order in the symmetric group $\mathrm{Sym}(\Lambda)$ on $\Lambda$.

It also holds that each $\rho \in \mathcal{A}_\mathrm{N}$ acts on $\mathcal{C}$ as an automorphism by $y \mapsto \rho(y)$ and $w_y^s \mapsto w_{\rho(y)}^s$.
In particular, we have $(\rho\rho')(x_I)=\rho(\rho'(x_I)) \in \mathcal{V}(C)$ for all $\rho,\rho' \in \mathcal{A}_\mathrm{N}$, where we put $\rho\rho'=\rho \circ \rho'$, therefore $\mathcal{A}_\mathrm{N}$ is a subgroup of $\mathrm{Sym}(\Lambda)$.
Thus this group $\mathcal{A}_\mathrm{N}$ is embedded into $\mathrm{Aut}\,\mathcal{C}$.
Moreover, for a $\rho \in \mathcal{A}_\mathrm{N}$ and any path $p$ in $\mathcal{C}$, the two paths $p$ and $\rho(p)$ represent the same element in $W$.

Now define the following maps
\begin{displaymath}
\mathsf{A}_\mathrm{N}:N_I \to \mathcal{A}_\mathrm{N},\ w \mapsto \rho \mbox{ for any } w \in C_{x_I,\rho(x_I)}
\end{displaymath}
and
\begin{displaymath}
h:\mathcal{A}_\mathrm{N} \to N_I,\ \rho \mapsto h_\rho=p_{x_I,\rho(x_I)} \enspace.
\end{displaymath}
Note that $\mathsf{A}_\mathrm{N}$ is a group homomorphism since
\begin{displaymath}
C_{x_I,\rho_1(x_I)}C_{x_I,\rho_2(x_I)}=C_{x_I,\rho_1(x_I)}C_{\rho_1(x_I),\rho_1\rho_2(x_I)} \subseteq C_{x_I,\rho_1\rho_2(x_I)} \mbox{ holds in } W \enspace.
\end{displaymath}
Then the following analogous results are deduced in a similar way:
\begin{thm}
\label{thm:exactsequence_N}
The sequence $1 \to C_I \hookrightarrow N_I \overset{\mathsf{A}_\mathrm{N}}{\to} \mathcal{A}_\mathrm{N} \to 1$ is exact and satisfies that $\mathsf{A}_\mathrm{N} \circ h=\mathrm{id}_{\mathcal{A}_\mathrm{N}}$.
If we define a subgroup $\widetilde{Y}_I$ of $N_I$ by
\begin{displaymath}
\widetilde{Y}_I=\{w \in N_I \mid (\Phi^{\perp I})^+=w \cdot (\Phi^{\perp I})^+\} \enspace,
\end{displaymath}
then $h(\mathcal{A}_\mathrm{N}) \subseteq \widetilde{Y}_I$, the sequence $1 \to Y_I \hookrightarrow \widetilde{Y}_I \overset{\mathsf{A}_\mathrm{N}}{\to} \mathcal{A}_\mathrm{N} \to 1$ is exact and $N_I=W^{\perp I} \rtimes \widetilde{Y}_I$, where $\widetilde{Y}_I$ acts on $W^{\perp I}$ as automorphisms of Coxeter system.
Hence
\begin{displaymath}
N_W(W_I)=W_I \rtimes (W^{\perp I} \rtimes \widetilde{Y}_I)=(W_I \times W^{\perp I}) \rtimes \widetilde{Y}_I \enspace.
\end{displaymath}
Moreover, if $q \in \pi_1(\mathcal{Y}^1;x_I)$, then we have
\begin{eqnarray}
\label{eq:exactsequence_N_1}
h_\rho^{-1}qh_\rho&=&\rho^{-1}(q_{(\rho(x_I))}) \enspace,\\
\nonumber
h_{\rho_1}h_{\rho_2}h_{\rho_3} \cdots h_{\rho_k}&=&p_{x_I,\rho_1(x_I)}\rho_{\left[1\right]}(p_{x_I,\rho_2(x_I)})\rho_{\left[2\right]}(p_{x_I,\rho_3(x_I)}) \cdots\\
\label{eq:exactsequence_N_2}
&&\qquad \quad \cdots \rho_{\left[k-1\right]}(p_{x_I,\rho_k(x_I)})p_{\rho_{\left[k\right]}(x_I),x_I}h_{\rho_{\left[k\right]}} \enspace,
\end{eqnarray}
where we put $\rho_{\left[i\right]}=\rho_1\rho_2 \cdots \rho_i$ for each $i$.
\end{thm}
\begin{prop}
\label{prop:Nsplits}
If the maximal tree $\mathcal{T}$ is stable under all the $\rho \in \mathcal{A}_\mathrm{N}$, then $h$ is a group homomorphism.
Hence in this case, the exact sequences of Theorem \ref{thm:exactsequence_N} split via the $h$.
\end{prop}
\begin{thm}
[See {\cite[Proposition 10.1]{Joh}}]
\label{thm:presentationofYtilde}
For the group $\mathcal{A}_{\mathrm{N}}$, choose a generating set $\mathcal{A}'_\mathrm{N}$ and fundamental relations of the form \lq\lq $\rho_1\rho_2 \cdots \rho_k=1$'' with $\rho_i \in \mathcal{A}'_\mathrm{N}$ (this is possible since each $\rho \in \mathcal{A}_\mathrm{N}$ has finite order).
Then $\widetilde{Y}_I$ admits a presentation, with generators given by $(w_\xi)_{(x_I)}$ for all $w_\xi \in E(\mathcal{Y})$ and $h_\rho$ for all $\rho \in \mathcal{A}'_\mathrm{N}$, and fundamental relations given by those of $Y_I$, relations (\ref{eq:exactsequence_N_1}) for $\rho \in \mathcal{A}'_\mathrm{N}$ and $q=(w_\xi)_{(x_I)}$ (where $w_\xi \in E(\mathcal{Y})$), and relations (\ref{eq:exactsequence_N_2}) for each sequence $\rho_1,\rho_2,\dots,\rho_k \in \mathcal{A}'_\mathrm{N}$ such that \lq\lq $\rho_1\rho_2 \cdots \rho_k=1$'' is one of the chosen fundamental relations of $\mathcal{A}_\mathrm{N}$.
\end{thm}
\begin{prop}
\label{prop:actiononWperp_Ytilde}
Let $w_y^s=s_{\gamma(y,s)}$ be a loop generator of $C$, and let $w,u \in Y_I$ and $\rho \in \mathcal{A}_\mathrm{N}$.
Then we have
\begin{displaymath}
ur(w,w_y^s)u^{-1}=r(uw,w_y^s) \mbox{ and } h_\rho r(w,w_y^s)h_\rho^{-1}=r(\rho(wp_{x_I,y})_{(x_I)},w_{\rho(y)}^s) \enspace.
\end{displaymath}
\end{prop}
%%
%%%%%
\section{Example}
\label{sec:example}
\begin{figure}[htb]
\centering
\begin{picture}(220,40)
\put(10,20){\circle*{10}}\put(10,3){\hbox to0pt{\hss$s_1$\hss}}
\put(50,20){\circle{10}}\put(50,3){\hbox to0pt{\hss$s_2$\hss}}
\put(90,20){\circle*{10}}\put(90,3){\hbox to0pt{\hss$s_3$\hss}}
\put(130,20){\circle*{10}}\put(130,3){\hbox to0pt{\hss$s_4$\hss}}
\put(170,20){\circle{10}}\put(170,3){\hbox to0pt{\hss$s_5$\hss}}
\put(210,20){\circle{10}}\put(210,3){\hbox to0pt{\hss$s_6$\hss}}
\put(15,20){\line(1,0){30}}
\put(55,20){\line(1,0){30}}\put(70,25){\hbox to0pt{\hss$4$\hss}}
\put(95,20){\line(1,0){30}}
\put(135,20){\line(1,0){30}}
\put(175,20){\line(1,0){30}}
\end{picture}
\caption{Coxeter graph for the example (here black circles signify elements of $I$)}
\label{fig:Coxetergraph}
\end{figure}
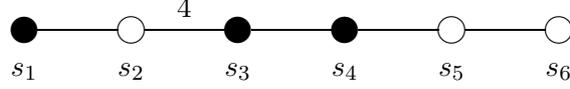

Let $(W,S)$ be the Coxeter system corresponding to the Coxeter graph in Figure \ref{fig:Coxetergraph}, and put $I=\{s_1,s_3,s_4\}$ and $x_I=(s_1,s_3,s_4)$.
Then we compute $Z_W(W_I)$ and $N_W(W_I)$.
%%%%%

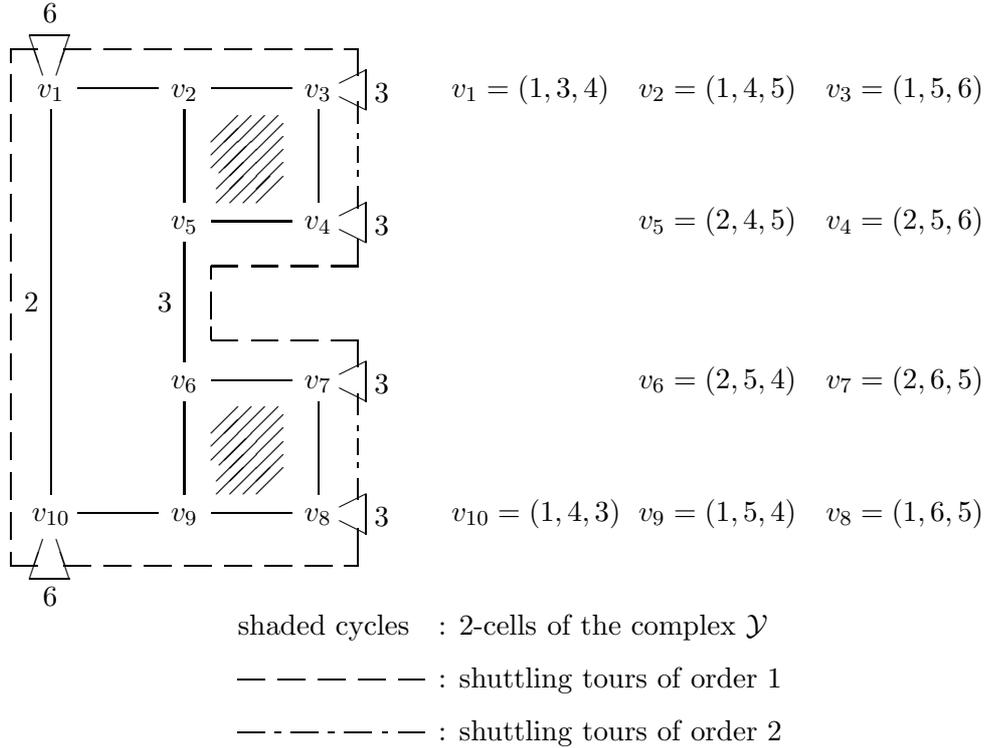
\begin{figure}[htbp]
\centering
\begin{picture}(340,230)
\put(10,195){\lower3pt\hbox to0pt{\hss$v_1$\hss}}
\put(60,195){\lower3pt\hbox to0pt{\hss$v_2$\hss}}
\put(110,195){\lower3pt\hbox to0pt{\hss$v_3$\hss}}
\put(110,145){\lower3pt\hbox to0pt{\hss$v_4$\hss}}
\put(60,145){\lower3pt\hbox to0pt{\hss$v_5$\hss}}
\put(60,85){\lower3pt\hbox to0pt{\hss$v_6$\hss}}
\put(110,85){\lower3pt\hbox to0pt{\hss$v_7$\hss}}
\put(110,35){\lower3pt\hbox to0pt{\hss$v_8$\hss}}
\put(60,35){\lower3pt\hbox to0pt{\hss$v_9$\hss}}
\put(10,35){\lower3pt\hbox to0pt{\hss$v_{10}$\hss}}
\put(20,195){\line(1,0){30}}
\put(70,195){\line(1,0){30}}
\put(70,145){\line(1,0){30}}
\put(70,85){\line(1,0){30}}
\put(20,35){\line(1,0){30}}
\put(70,35){\line(1,0){30}}
\put(10,187){\line(0,-1){145}}\put(0,111){$2$}
\put(60,187){\line(0,-1){35}}
\put(110,187){\line(0,-1){35}}
\put(60,137){\line(0,-1){45}}\put(50,111){$3$}
\put(60,77){\line(0,-1){35}}
\put(110,77){\line(0,-1){35}}
\put(7,200){\line(-1,3){5}}\put(12,200){\line(1,3){5}}\put(2,215){\line(1,0){15}}\put(7,220){$6$}
\put(7,25){\line(-1,-3){5}}\put(12,25){\line(1,-3){5}}\put(2,10){\line(1,0){15}}\put(7,0){$6$}
\put(118,197){\line(2,1){10}}\put(118,192){\line(2,-1){10}}\put(128,202){\line(0,-1){15}}\put(131,190){$3$}
\put(118,147){\line(2,1){10}}\put(118,142){\line(2,-1){10}}\put(128,152){\line(0,-1){15}}\put(131,140){$3$}
\put(118,87){\line(2,1){10}}\put(118,82){\line(2,-1){10}}\put(128,92){\line(0,-1){15}}\put(131,80){$3$}
\put(118,37){\line(2,1){10}}\put(118,32){\line(2,-1){10}}\put(128,42){\line(0,-1){15}}\put(131,30){$3$}
\put(80,185){\line(-1,-1){10}}
\put(85,185){\line(-1,-1){15}}
\put(90,185){\line(-1,-1){20}}
\put(95,185){\line(-1,-1){22}}
\put(97,182){\line(-1,-1){25}}
\put(97,177){\line(-1,-1){25}}
\put(97,172){\line(-1,-1){20}}
\put(97,167){\line(-1,-1){15}}
\put(97,162){\line(-1,-1){10}}
\put(80,75){\line(-1,-1){10}}
\put(85,75){\line(-1,-1){15}}
\put(90,75){\line(-1,-1){20}}
\put(95,75){\line(-1,-1){22}}
\put(97,72){\line(-1,-1){25}}
\put(97,67){\line(-1,-1){25}}
\put(97,62){\line(-1,-1){20}}
\put(97,57){\line(-1,-1){15}}
\put(97,52){\line(-1,-1){10}}
\put(5,210){\line(-1,0){10}}\put(-5,210){\line(0,-1){5}}\put(-5,200){\line(0,-1){10}}\put(-5,185){\line(0,-1){10}}\put(-5,170){\line(0,-1){10}}\put(-5,155){\line(0,-1){10}}\put(-5,140){\line(0,-1){10}}\put(-5,125){\line(0,-1){10}}\put(-5,110){\line(0,-1){10}}\put(-5,95){\line(0,-1){10}}\put(-5,80){\line(0,-1){10}}\put(-5,65){\line(0,-1){10}}\put(-5,50){\line(0,-1){10}}\put(-5,35){\line(0,-1){10}}\put(-5,20){\line(0,-1){5}}\put(-5,15){\line(1,0){10}}
\put(15,210){\line(1,0){10}}\put(30,210){\line(1,0){10}}\put(45,210){\line(1,0){10}}\put(60,210){\line(1,0){10}}\put(75,210){\line(1,0){10}}\put(90,210){\line(1,0){10}}\put(105,210){\line(1,0){10}}\put(120,210){\line(1,0){5}}\put(125,210){\line(0,-1){10}}
\put(15,15){\line(1,0){10}}\put(30,15){\line(1,0){10}}\put(45,15){\line(1,0){10}}\put(60,15){\line(1,0){10}}\put(75,15){\line(1,0){10}}\put(90,15){\line(1,0){10}}\put(105,15){\line(1,0){10}}\put(120,15){\line(1,0){5}}\put(125,15){\line(0,1){14}}
\put(125,138){\line(0,-1){10}}\put(125,128){\line(-1,0){10}}\put(110,128){\line(-1,0){10}}\put(95,128){\line(-1,0){10}}\put(80,128){\line(-1,0){10}}\put(70,128){\line(0,-1){5}}\put(125,90){\line(0,1){10}}\put(125,100){\line(-1,0){10}}\put(110,100){\line(-1,0){10}}\put(95,100){\line(-1,0){10}}\put(80,100){\line(-1,0){10}}\put(70,100){\line(0,1){5}}\put(70,118){\line(0,-1){10}}
\put(125,190){\line(0,-1){8}}\put(125,178){\line(0,-1){2}}\put(125,173){\line(0,-1){6}}\put(125,162){\line(0,1){2}}\put(125,150){\line(0,1){8}}
\put(125,80){\line(0,-1){8}}\put(125,68){\line(0,-1){2}}\put(125,63){\line(0,-1){6}}\put(125,52){\line(0,1){2}}\put(125,40){\line(0,1){8}}
\put(160,192){$v_1=(1,3,4)$}\put(230,192){$v_2=(1,4,5)$}\put(300,192){$v_3=(1,5,6)$}
\put(230,142){$v_5=(2,4,5)$}\put(300,142){$v_4=(2,5,6)$}
\put(230,82){$v_6=(2,5,4)$}\put(300,82){$v_7=(2,6,5)$}
\put(160,32){$v_{10}=(1,4,3)$}\put(230,32){$v_9=(1,5,4)$}\put(300,32){$v_8=(1,6,5)$}
\end{picture}
\begin{picture}(180,50)
\put(0,40){shaded cycles}\put(75,40){: $2$-cells of the complex $\mathcal{Y}$}
\put(0,23){\line(1,0){10}}\put(15,23){\line(1,0){10}}\put(30,23){\line(1,0){10}}\put(45,23){\line(1,0){10}}\put(60,23){\line(1,0){10}}\put(75,20){: shuttling tours of order $1$}
\put(0,3){\line(1,0){10}}\put(14,3){\line(1,0){2}}\put(20,3){\line(1,0){10}}\put(34,3){\line(1,0){2}}\put(40,3){\line(1,0){10}}\put(54,3){\line(1,0){2}}\put(60,3){\line(1,0){10}}\put(75,0){: shuttling tours of order $2$}
\end{picture}
\caption{Graph $\mathcal{C}$, complex $\mathcal{Y}$ and shuttling tours for the example}
\label{fig:complex}
\end{figure}

First we construct the graph $\mathcal{C}$.
We start with a vertex $x=x_I$, and choose $s \in S \smallsetminus [x]$ such that $[x]_{\sim s}$ is of finite type.
Say, $s=s_5$, hence $[x]_{\sim s}=\{s_3,s_4,s_5\}$ and $[x]_{\sim s} \smallsetminus \{s\}=\{s_3,s_4\}$ are of type $A_3$ and $A_2$, respectively; therefore
\begin{displaymath}
\alpha_{s_3} \overset{w_0([x]_{\sim s} \smallsetminus \{s\})}{\mapsto} -\alpha_{s_4} \overset{w_0([x]_{\sim s})}{\mapsto} \alpha_{s_4} \mbox{ and } 
\alpha_{s_4} \overset{w_0([x]_{\sim s} \smallsetminus \{s\})}{\mapsto} -\alpha_{s_3} \overset{w_0([x]_{\sim s})}{\mapsto} \alpha_{s_5} \enspace.
\end{displaymath}
Thus we have $\varphi(x,s)=((s_1,s_4,s_5),s_3)$, hence we add to the graph $\mathcal{C}$ a new vertex $y=(s_1,s_4,s_5)$ and a new edge $w_x^s$ from $x$ to $y$.
Iterating such a process, we obtain the entire graph $\mathcal{C}$ finally as depicted in Figure \ref{fig:complex} by solid lines.
In the figure, we abbreviate $s_i$ to $i$, and each edge $w_y^s$ is labelled $s$ when the tuples $w_y^s \cdot y$ and $y$ consist of the same contents, since now the element $s$ cannot be determined by the vertices $y$ and $w_y^s \cdot y$ only.
%%%%%
%%%%%

Concerning the circular tours and the shuttling tours, choose a subset $J \subseteq S$ with $|J|=|I|+2=5$, say $J=S \smallsetminus \{s_3\}$, and a connected component $\mathcal{G}$ of $\mathcal{C}(J)$, say the cycle $v_2v_3v_4v_5v_2$, such that $J_{\sim (J \smallsetminus [y])}$ is of finite type for a vertex $y$ of $\mathcal{G}$.
Since $\mathcal{G}$ is now a cycle, it contains a circular tour by definition, yielding a $2$-cell of $\mathcal{Y}$ with boundary $\mathcal{G}$.
On the other hand, if $J=S \smallsetminus \{s_2\}$, and $\mathcal{G}$ consists of two loops $w_{v_1}^{s_6}$ and $w_{v_3}^{s_3}$ and a path $p=v_1v_2v_3$, then $J_{\sim (J \smallsetminus [v_1])}=J \smallsetminus \{s_1\}$ is of finite type, hence we obtain a shuttling tour $w_{v_1}^{s_6}pw_{v_3}^{s_3}p^{-1}$ of order one (see the first row of Table \ref{tab:listofedges_1}, where $K=J \smallsetminus \{s_1\}$).
This $p$ is depicted in Figure \ref{fig:complex} as a broken line.
Iteration of such a process enumerates the circular tours (hence the $2$-cells of $\mathcal{Y}$) and the shuttling tours.
(Note that the cycle $v_1v_2v_5v_6v_9v_{10}$ is \emph{not} the boundary of a $2$-cell, since $J_{\sim (J \smallsetminus \{s_1\})}=J$ is \emph{not} of finite type where $J=S \smallsetminus \{s_6\}$.)

We determine the structure of $B_I \subseteq Z_W(W_I)$ and $\widetilde{Y}_I \subseteq N_W(W_I)$.
Let $e_1$, $e_2$, $e_3$ be the edges $v_4 \leftarrow v_3$, $v_6 \leftarrow v_5$, $v_8 \leftarrow v_7$, respectively, and take a unique maximal tree $\mathcal{T}$ in $\mathcal{Y}^1$ such that $e_1,e_2,e_3 \not\in \mathcal{T}$.
By Theorems \ref{thm:presentationofY} and \ref{thm:presentationofpi1Y}, $Y_I \simeq \pi_1(\mathcal{Y};x_I)$ is generated by the three $(e_i)_{(v_1)}$, and has fundamental relations \lq\lq $(e_1)_{(v_1)}=1$'' and \lq\lq $(e_3)_{(v_1)}=1$'' induced by the two $2$-cells of $\mathcal{Y}$.
As a result, $Y_I$ is a free group of rank one (i.e., an infinite cyclic group) generated by $a=(e_2)_{(v_1)}$.
Now for $B_I$, we have $x_I{}^A=v_{10} \in \mathcal{V}(C)$ and $\mathcal{A}=\{\emptyset,A\}$, where $A=\{2,3\} \subseteq \Lambda=\{1,2,3\}$.
Now $\mathcal{T}$ is stable under the map $\tau_A$ that sends each $v_i$ to $v_{11-i}$, therefore Proposition \ref{prop:Zsplits} implies that $B_I$ is a semidirect product of $Y_I$ by a subgroup $\langle g_A \rangle \simeq \mathcal{A} \simeq \{\pm 1\}$.
Moreover, by the formulae in Theorem \ref{thm:exactsequence_Z}, we have
\begin{displaymath}
g_Aag_A^{-1}=\tau_A(a)_{(v_1)}={e_2^{-1}}_{(v_1)}=a^{-1} \enspace.
\end{displaymath}
Hence, by putting $b=b'=g_A$ and $a'=ba$, we have
\begin{displaymath}
B_I=\langle a \rangle \rtimes \langle b \rangle \simeq \mathbb{Z} \rtimes \{\pm 1\} \mbox{ and } B_I=\langle a',b' \mid a'{}^2=b'{}^2=1 \rangle \simeq W(\widetilde{A_1}) \enspace,
\end{displaymath}
an affine Coxeter group of type $\widetilde{A_1}$.

Similarly, putting $\rho=\genfrac{(}{)}{0pt}{}{1\,2\,3}{1\,3\,2} \in \mathrm{Sym}(\Lambda)$, we have $\mathcal{A}_\mathrm{N}=\{\mathrm{id}_\Lambda,\rho\}$, $\rho(v_i)=v_{11-i}$ and $\mathcal{T}$ is stable under the $\rho$.
Thus Theorem \ref{thm:exactsequence_N} and Proposition \ref{prop:Nsplits} imply that
\begin{displaymath}
h_\rho^{-1}ah_\rho=\rho^{-1}((e_2)_{(v_{10})})=(e_2)_{(v_1)}{}^{-1}=a^{-1}
\end{displaymath}
and we have $\widetilde{Y}_I=\langle a \rangle \rtimes \langle h_\rho \rangle \simeq W(\widetilde{A_1})$.

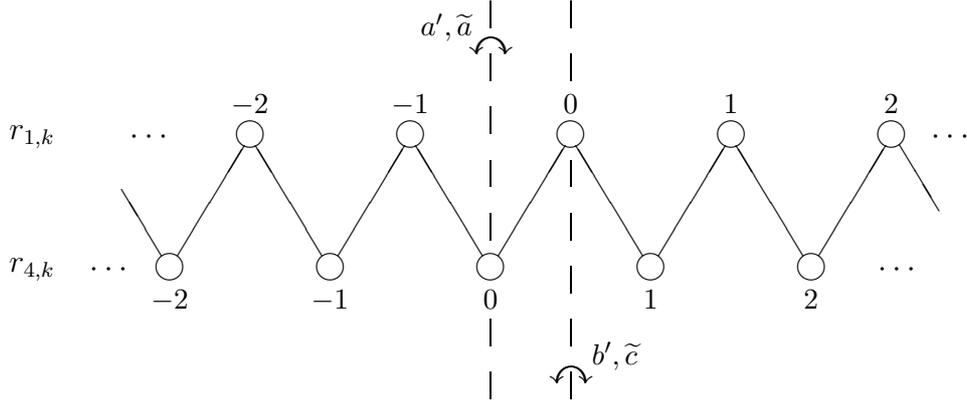
\begin{figure}[htb]
\centering
\begin{picture}(360,160)(10,0)
\put(100,100){\circle{10}}\put(100,108){\hbox to0pt{\hss$-2$\hss}}
\put(160,100){\circle{10}}\put(160,108){\hbox to0pt{\hss$-1$\hss}}
\put(220,100){\circle{10}}\put(220,108){\hbox to0pt{\hss$0$\hss}}
\put(280,100){\circle{10}}\put(280,108){\hbox to0pt{\hss$1$\hss}}
\put(340,100){\circle{10}}\put(340,108){\hbox to0pt{\hss$2$\hss}}
\put(70,50){\circle{10}}\put(70,34){\hbox to0pt{\hss$-2$\hss}}
\put(130,50){\circle{10}}\put(130,34){\hbox to0pt{\hss$-1$\hss}}
\put(190,50){\circle{10}}\put(190,34){\hbox to0pt{\hss$0$\hss}}
\put(250,50){\circle{10}}\put(250,34){\hbox to0pt{\hss$1$\hss}}
\put(310,50){\circle{10}}\put(310,34){\hbox to0pt{\hss$2$\hss}}
\put(67,54){\line(-3,5){15}}
\put(73,54){\line(3,5){25}}
\put(127,54){\line(-3,5){25}}
\put(133,54){\line(3,5){25}}
\put(187,54){\line(-3,5){25}}
\put(193,54){\line(3,5){25}}
\put(247,54){\line(-3,5){25}}
\put(253,54){\line(3,5){25}}
\put(307,54){\line(-3,5){25}}
\put(313,54){\line(3,5){25}}
\put(343,96){\line(3,-5){15}}
\put(10,98){\large $r_{1,k}$}\put(55,96){\large $\cdots$}\put(355,96){\large $\cdots$}
\put(10,48){\large $r_{4,k}$}\put(40,46){\large $\cdots$}\put(335,46){\large $\cdots$}
\put(220,0){\line(0,1){10}}\put(220,20){\line(0,1){10}}\put(220,40){\line(0,1){10}}\put(220,60){\line(0,1){10}}\put(220,80){\line(0,1){10}}\put(220,120){\line(0,1){10}}\put(220,140){\line(0,1){10}}
\put(190,0){\line(0,1){10}}\put(190,20){\line(0,1){10}}\put(190,60){\line(0,1){10}}\put(190,80){\line(0,1){10}}\put(190,100){\line(0,1){10}}\put(190,120){\line(0,1){10}}\put(190,140){\line(0,1){10}}
\put(214,6){\Large $\curvearrowright$}\put(212,6){\Large $\curvearrowleft$}\put(228,13){$b',\widetilde{c}$}
\put(184,130){\Large $\curvearrowright$}\put(182,130){\Large $\curvearrowleft$}\put(164,137){$a',\widetilde{a}$}
\end{picture}
\caption{Structure of the centralizer in the example (here two adjacent generators commute and the product of two non-adjacent generators has infinite order)}
\label{fig:centralizer}
\end{figure}

We compute the structure of $W^{\perp I}$ by using Theorems \ref{thm:presentationofWperp} and \ref{thm:relationinWperp}.
Let $\xi_i$ (with $i=1,3,4,7,8,10$) denote the unique pair of the form $(v_i,s_j)$ such that $w_{\xi_i}$ is a loop in $\mathcal{C}$.
Then $W^{\perp I}$ is generated by the $r(a^k,w_{\xi_i})$ for all $i$ and $k \in \mathbb{Z}$.
For their relations, let $w_{\xi_i}qw_{\xi_j}q^{-1}$ be one of the six shuttling tours.
If $(i,j)=(1,10)$, $(1,3)$ or $(8,10)$, then the shuttling tour has order one and $q \subseteq \mathcal{T}$, therefore we have
\begin{displaymath}
(a^k,w_{\xi_i}) \overset{1}{\sim} (a^kq_{(v_1)},w_{\xi_j})=(a^k,w_{\xi_j}) \enspace.
\end{displaymath}
If $(i,j)=(4,7)$, then the order is also one and $q_{(v_1)}=a^{-1}$, therefore we have
\begin{displaymath}
(a^k,w_{\xi_4}) \overset{1}{\sim} (a^kq_{(v_1)},w_{\xi_7})=(a^{k-1},w_{\xi_7}) \enspace.
\end{displaymath}
Moreover, if $(i,j)=(3,4)$ or $(7,8)$, then the order is two and $q \subseteq \mathcal{T}$, therefore we have
\begin{displaymath}
(a^k,w_{\xi_i}) \overset{2}{\sim} (a^kq_{(v_1)},w_{\xi_j})=(a^k,w_{\xi_j}) \enspace.
\end{displaymath}
These exhaust the relations $\overset{m}{\sim}$ on the set $\mathcal{R}_{x_I}$.
Thus the generating set $R^I$ of $W^{\perp I}$ consists of infinitely many \emph{distinct} elements
\begin{displaymath}
r_{1,k}=r(a^k,w_{\xi_1}) \mbox{ and } r_{4,k}=r(a^k,w_{\xi_4})\,,\, \mbox{ with } k \in \mathbb{Z}
\end{displaymath}
and the fundamental relations are as depicted in Figure \ref{fig:centralizer} where, differently from usual Coxeter graphs, any two adjacent generators commute and a product of any two non-adjacent generators has infinite order.
In particular, $W^{\perp I}$ is \emph{not} finitely generated though $W$ is finitely generated.
Note that this $W^{\perp I}$ possesses no finite irreducible component; this fact can also be deduced directly by using Theorem \ref{thm:conditionfornonfinitepart}.
%%%%%
%%%%%

Finally, we determine the actions of $B_I$ and $\widetilde{Y}_I$ on $W^{\perp I}$ by using Propositions \ref{prop:actiononWperp_B} and \ref{prop:actiononWperp_Ytilde}.
Put $c=h_\rho$, $\widetilde{a}=ca$ and $\widetilde{c}=c$.
Note that $\{a',b'\}$ and $\{\widetilde{a},\widetilde{c}\}$ are the generating sets of $B_I$ and $\widetilde{Y}_I$, respectively, as affine Coxeter groups of type $\widetilde{A_1}$.
Then by the formulae given in those propositions, we have
\begin{eqnarray*}
a'r_{1,k}a'{}^{-1}=\widetilde{a}r_{1,k}\widetilde{a}^{-1}=r_{1,-k-1}&,&
a'r_{4,k}a'{}^{-1}=\widetilde{a}r_{4,k}\widetilde{a}^{-1}=r_{4,-k} \enspace,\\
b'r_{1,k}b'{}^{-1}=\widetilde{c}r_{1,k}\widetilde{c}^{-1}=r_{1,-k}&,&
b'r_{4,k}b'{}^{-1}=\widetilde{c}r_{4,k}\widetilde{c}^{-1}=r_{4,1-k} \enspace.
\end{eqnarray*}
Thus $a'$ and $\widetilde{a}$ act on the graph in Figure \ref{fig:centralizer} as the rotation round the vertical axis through $r_{4,0}$, while $b'$ and $\widetilde{c}$ do as the rotation round the vertical axis through $r_{1,0}$.
This determines the actions of $B_I$ and $\widetilde{Y}_I$.

Finally, we have $Z(W_I)=\langle s_1 \rangle \simeq \{\pm 1\}$ (see Section \ref{sec:longestelement}).
Hence the structure of $Z_W(W_I)$ and $N_W(W_I)$ have been completely determined.

\noindent
\textbf{Koji Nuida}\\
Present address: Research Center for Information Security (RCIS), National Institute of Advanced Industrial Science and Technology (AIST), AIST Tsukuba Central 2, 1-1-1 Umezono, Tsukuba, Ibaraki 305-8568, Japan\\
E-mail: k.nuida[at]aist.go.jp
\end{document}